# Mordukhovich derivatives of the metric projection operator in Hilbert spaces

Jinlu Li

Department of Mathematics
Shawnee State University
Portsmouth, Ohio 45662 USA
Email: jli@shawnee.edu

**Abstract**

In this paper, we study the generalized differentiability of the metric projection operator in Hilbert spaces. We find exact expressions for Mordukhovich derivatives for the metric projection operator onto closed balls in Hilbert spaces and positive cones in Euclidean spaces and in real Hilbert space $l_2$. We investigate the connections between Fréchet differentiability, Gâteaux directional differentiability and the Mordukhovich derivatives of the metric projection operator in Hilbert spaces.

1.  Introduction

Throughout this paper, let $(H, \|\cdot\|)$ be a real Hilbert space with inner product $\langle \cdot, \cdot \rangle$ and with the origin $\theta$. Let $f: H \to H$ be a single-valued mapping. The continuity of $f$ is an important topic in the analysis in Hilbert spaces. In addition to continuity, the smoothness of single-valued mappings in Hilbert spaces has been studied with respect to several types of differentiability, such as, the Gâteaux directional differentiability, Fréchet differentiability, and strict Fréchet differentiability.

It is well-known that set-valued mappings play very important roles in nonlinear analysis. To extend the differentiability from single-valued mappings to set-valued mappings, the generalized differentiability of set-valued mappings has been introduced (see Mordukhovich [11]). The generalized differentiability of set-valued mappings has been applied to optimization theory, approximation theory, theory of control, so and forth (see [11]).

In this paper, we study the Mordukhovich derivatives of the metric projection operator in Hilbert spaces, which is considered as a special case of set-valued mappings. More precisely speaking, let $C$ be a nonempty closed and convex subset of $H$ and $P_C: H \to C$ denote the metric projection operator, which is defined, for any $x \in H$, by

$$\|x - P_C(x)\| \leq \|x - z\|, \text{ for all } z \in C.$$

$P_C$ is a well-defined single-valued mapping that has many useful properties, such as continuity, monotonicity, and non-expansiveness. These properties have been studied by many authors and have been applied to many fields in analysis in Hilbert spaces (see [1]).

The smoothness of the metric projection operator $P_C$ in Hilbert spaces has been described by several types of differentiability, such as, (Gâteaux) directional differentiability, Fréchet differentiability and strict Fréchet differentiability (see [1, 3−5, 7−8, 12−14]). It is worth to note that the differentiability of the metric projection operator has been extended from Hilbert spaces to uniformly convex and uniformly smooth Banach spaces (see [2, 6]). The differentiability of the metric projection operator has been applied to approximation theory, optimization theory, and variational inequalities in Hilbert spaces (see [1, 9−10]).

In section 2, we recall some concepts and properties of (Gâteaux) directional differentiability, Fréchet differentiability, strict Fréchet differentiability and generalized differentiability of single-valued mappings in Hilbert spaces. In section 3, we investigate the Mordukhovich derivatives of the metric projection operator onto closed balls with center origin. In section 4, we investigate the Mordukhovich derivatives of the metric projection operator onto the positive cones in Euclidean spaces. In section 5, we consider the real Hilbert space $l_2$.

For a single-valued mapping $f: H \to H$, if $f$ is Fréchet differentiable at some point $x$ in $H$, by definition, the Fréchet derivative of $f$ at this point is a linear and continuous mapping from $H$ to $H$. If $f$ is Gâteaux directionally differentiable at some point $x$ along some direction in $H$, by definition, the Gâteaux directional derivative of $f$ at this point with this direction is a single point in $H$. However, for a single-valued mapping $f$, the Mordukhovich derivatives of $f$ at $(x, f(x))$ is a set. In part (iii) of Theorem 5.1 in section 5, we show that the Mordukhovich derivatives of the metric projection onto the positive cone in $l_2$ is a set that is not a singleton.

## 2. Preliminaries

### 2.1 Fréchet differentiability of single-valued mappings in Hilbert spaces

In this subsection, we recall the definitions and some properties of Fréchet and generalized differentiability of single-valued mappings in Hilbert spaces. In this paper, we will closely adapt the notations from [1] and [11]. However, since we consider Hilbert spaces in this paper instead of general Banach spaces and we only study single-valued mappings, so, we will slightly simplify the notations and concepts used in [1] and [11].

We first recall the concepts of the differentiability of single-valued mappings from $H$ to $H$.

**Definition 3.1 in [6]**. (Gâteaux directional differentiability of single-valued mappings in Hilbert spaces). Let $f: H \to H$ be a single-valued mapping. For $x \in H$ and $w \in H$ with $w \neq \theta$, if the following limit exists, which is a point in $H$,

$$f'(x)(w) = \lim_{t \downarrow 0} \frac{f(x+tw) - f(x)}{t}, \tag{2.1}$$

then, $f$ is said to be Gâteaux directionally differentiable at point $x$ along direction $w$. $f'(x)(w)$ is called the (Gâteaux) directional derivative of $f$ at point $x$ along direction $w$. Let $A$ be a subset in

$H$. If $f$ is (Gâteaux) directionally differentiable at every point $x \in A$, then $f$ is said to be (Gâteaux) directionally differentiable on $A \subseteq H$.

**Definition 1.13 in [11]**. (Fréchet derivatives of single-valued mappings in Hilbert spaces). Let $f: H \to H$ be a single-valued mapping and let $\bar{x} \in H$. If there is a linear and continuous mapping $\nabla f(\bar{x}): H \to H$ such that

(i) $\quad \lim_{x \to \bar{x}} \dfrac{f(x) - f(\bar{x}) - \nabla f(\bar{x})(x - \bar{x})}{\|x - \bar{x}\|} = \theta.$

then $f$ is said to be Fréchet differentiable at $\bar{x}$ and $\nabla f(\bar{x})$ is called the Fréchet derivative of $f$ at $\bar{x}$;

(ii) $\quad$ If $\lim_{u \to \bar{x}, v \to \bar{x}} \dfrac{f(u) - f(v) - \nabla f(\bar{x})(u - v)}{\|u - v\|} = \theta,$

then $f$ is said to be strictly Fréchet differentiable at $\bar{x}$.

### 2.2 Mordukhovich derivatives of single-valued mappings in Hilbert spaces

In [11], the concepts of generalized differentiation of set-valued mappings in Banach spaces are introduced. Here, we rewrite those definitions with respect to single-valued mappings in Hilbert spaces.

**Definition 1.1 in [11]** (generalized normals in Hilbert spaces). Let $\Omega$ be a nonempty subset of a Hilbert space $H$. For any $x \in H$, we define

(i) $\quad \widehat{N}(x; \Omega) = \left\{ z \in H : \limsup_{\substack{u \to x \\ \Omega}} \dfrac{\langle z, u - x \rangle}{\|u - x\|} \leq 0 \right\}$, for any $x \in \Omega$;

(ii) $\quad \widehat{N}(x; \Omega) = \emptyset$, for any $x \in H \setminus \Omega$.

The elements of $\widehat{N}(x; X)$ are called Fréchet normals and $\widehat{N}(x; X)$ is the prenormal cone to $\Omega$ at $x$. Then, based on the above definition, we recall the Mordukhovich derivatives of single-valued mappings in Hilbert spaces. Note that, in [11], Mordukhovich derivatives are named by precoderivatives or Fréchet coderivatives.

**Definition 1.32 in [11]** (Mordukhovich derivatives of single-valued mappings in Hilbert spaces). Let $f: H \to H$ be a single-valued mapping. The Mordukhovich derivatives of $f$ at $(x, f(x))$ is

$$\widehat{D}^* f(x, f(x))(y) := \widehat{D}^* f(x)(y) = \{z \in H : (z, -y) \in \widehat{N}((x, f(x)); \mathrm{gph} f)\}.$$

It follows that

$$\widehat{D}^* f(x, u)(y) = \emptyset, \text{ for all } y \in H \text{ if } u \neq f(x).$$

Next theorem provides the connections between Mordukhovich derivatives and Fréchet derivatives for single-valued mappings in Hilbert spaces.

**Theorem 1.38 in [11]**. *Let f: $H \to H$ be Fréchet differentiable at $\bar{x}$. Then, Mordukhovich derivatives satisfy the following equation*

$$\widehat{D}^*f(\bar{x})(y) = \{(\nabla f(\bar{x}))y\}, \text{ for all } y \in H.$$

Let $H \times H$ be the orthogonal product of $H$ equipped with the inner product $\langle \cdot, \cdot \rangle_{H \times H}$, which is abbreviated as $\langle \cdot, \cdot \rangle$. It is defined by,

$$\langle (x,y), (u,v) \rangle := \langle x, u \rangle + \langle y, v \rangle, \text{ for any } (x,y), (u,v) \in H \times H.$$

Let $\|\cdot\|_{H \times H}$ be the norm in $H \times H$ induced by the inner product $\langle \cdot, \cdot \rangle$. Then, we have

$$\|(x,y)\|_{H \times H}^2 = \langle (x,y), (x,y) \rangle = \langle x, x \rangle + \langle y, y \rangle = \|x\|^2 + \|y\|^2.$$

It follows that $(H \times H, \|\cdot\|_{H \times H})$ is a real Hilbert space with inner product $\langle \cdot, \cdot \rangle_{H \times H}$ and with the origin $(\theta, \theta)$, which is simply denoted by $\theta$.

For any single-valued mapping $f: H \to H$. By Definitions 1.1 and 1.32 in [11], for any $\bar{x} \in H$, the Mordukhovich derivatives of $f$ at $(\bar{x}, f(\bar{x}))$ is

$$\widehat{D}^*f(\bar{x}, f(\bar{x}))(y)$$

$$:= \widehat{D}^*f(\bar{x})(y)$$

$$= \{z \in H : (z, -y) \in \widehat{N}((\bar{x}, f(\bar{x})); \text{gph} f)\}$$

$$= \left\{ z \in H : \limsup_{(u,f(u)) \xrightarrow{\text{gph} f} (\bar{x}, f(\bar{x}))} \frac{\langle (z, -y), (u, f(u)) - (\bar{x}, f(\bar{x})) \rangle}{\|(u, f(u)) - (\bar{x}, f(\bar{x}))\|_{H \times H}} \leq 0 \right\}$$

$$= \left\{ z \in H : \limsup_{(u,f(u)) \xrightarrow{\text{gph} f} (\bar{x}, f(\bar{x}))} \frac{\langle (z, -y), (u, f(u)) - (\bar{x}, f(\bar{x})) \rangle}{\sqrt{\|u - \bar{x}\|^2 + \|f(u) - f(\bar{x})\|^2}} \leq 0 \right\}, \text{ for } y \in H.$$

The following lemma provides a useful formula to calculate the Fréchet coderivatives of single-valued mappings on Hilbert spaces.

**Lemma 2.1.** *Let $H$ be a Hilbert space. The Mordukhovich derivatives of $f$ at a point $\bar{x} \in H$ satisfies that, for any $y \in H$,*

$$\widehat{D}^*f(\bar{x})(y) = \left\{ z \in H : \limsup_{(u,f(u)) \xrightarrow{\text{gph} f} (\bar{x}, f(\bar{x}))} \frac{\langle z, u - \bar{x} \rangle - \langle y, f(u) - f(\bar{x}) \rangle}{\|u - \bar{x}\| + \|f(u) - f(\bar{x})\|} \leq 0 \right\}. \quad (2.1)$$

*In particular, if $f$ is continuous, then*

$$\widehat{D}^*f(\bar{x})(y) = \left\{ z \in H : \limsup_{u \to \bar{x}} \frac{\langle z, u - \bar{x} \rangle - \langle y, f(u) - f(\bar{x}) \rangle}{\|u - \bar{x}\| + \|f(u) - f(\bar{x})\|} \leq 0 \right\}. \quad (2.2)$$

*Proof.* By the definition of the inner product in $H \times H$, we have

$$\langle (z, -y),\ (u, f(u)) - (\bar{x}, f(\bar{x})) \rangle = \langle z, u - \bar{x} \rangle - \langle y,\ f(u) - f(\bar{x}) \rangle.$$

Notice that

$$\frac{\sqrt{2}}{2} \leq \frac{\sqrt{\|u-\bar{x}\|^2 + \|f(u) - f(\bar{x})\|^2}}{\|u-\bar{x}\| + \|f(u) - f(\bar{x})\|} \leq 1.$$

Then, by the definition of $\widehat{D}^* f(\bar{x})(y)$, this reduces (2.1) immediately. □

## 3. Mordukhovich derivatives of the metric projection onto balls in Hilbert spaces

### 3.1. An orthogonal decomposition in Hilbert spaces

For any $x, \bar{x} \in H$, as usual, we write $x \perp \bar{x}$ if and only if $\langle x, \bar{x} \rangle = 0$. For any $\bar{x} \in H \setminus \{\theta\}$, let $S(\bar{x})$ denote the one-dimensional subspace of $H$ generated by $\bar{x}$. Let $O(\bar{x})$ denote the orthogonal subspace of $\bar{x}$ (or $S(\bar{x})$) in $H$. $H$ has the following orthogonal decomposition

$$H = S(\bar{x}) \oplus O(\bar{x}).$$

More precisely speaking, for this given $\bar{x} \in H \setminus \{\theta\}$ and for any $x \in H$, $x$ enjoys the following orthogonal representation

$$x = \frac{\langle x, \bar{x} \rangle}{\|\bar{x}\|^2} \bar{x} + \left( x - \frac{\langle x, \bar{x} \rangle}{\|\bar{x}\|^2} \bar{x} \right), \text{ for all } x \in H.$$

By using the above orthogonal representations of elements in $H$, for this fixed $\bar{x} \in H \setminus \{\theta\}$, we define a real valued function $a(\bar{x}; \cdot) \colon H \to \mathbb{R}$ by

$$a(\bar{x}; x) := \frac{\langle x, \bar{x} \rangle}{\|\bar{x}\|^2}, \text{ for all } x \in H.$$

And a mapping $o(\bar{x}; \cdot) \colon H \to O(\bar{x})$ by

$$o(\bar{x}; x) := x - \frac{\langle x, \bar{x} \rangle}{\|\bar{x}\|^2} \bar{x}, \text{ for all } x \in H.$$

The following lemma provides some properties of $a(\bar{x}; \cdot)$ and $o(\bar{x}; \cdot)$. These properties will be repeatedly used in this section.

**Lemma 3.1 in [11]**. *For any given fixed $\bar{x} \in H \setminus \{\theta\}$, the real valued function $a(\bar{x}; \cdot)$ and the mapping $o(\bar{x}; \cdot)$ have the following properties.*

(i)  $a(\bar{x}; \cdot) \colon H \to \mathbb{R}$ *is a real valued linear and continuous function;*
(ii) $o(\bar{x}; \cdot) \colon H \to O(\bar{x})$ *is a linear and continuous mapping;*
(iii) $a(\bar{x}; \cdot)\bar{x}$ *and $o(\bar{x}; \cdot)$ are orthogonal of each other and for any $u, v \in H$, we have*

(a) $\langle u, v \rangle = a(\bar{x}; u) a(\bar{x}; v) \|\bar{x}\|^2 + \langle o(\bar{x}; u), o(\bar{x}; v) \rangle$;
(b) $\langle a(\bar{x}; u)\bar{x}, o(\bar{x}; v) \rangle = 0$;
(c) $\|a(\bar{x}; u)\bar{x} + o(\bar{x}; v)\|^2 = (a(\bar{x}; u))^2 \|\bar{x}\|^2 + \|o(\bar{x}; v)\|^2$;

(d) $\|u\|^2 = (a(\bar{x}; u))^2 \|\bar{x}\|^2 + \|o(\bar{x}; u)\|^2$;

(e) $\|u + v\|^2 = (a(\bar{x}; u) + a(\bar{x}; u))^2 \|\bar{x}\|^2 + \|o(\bar{x}; u) + o(\bar{x}; v)\|^2$.

(iv) $\quad u \to \bar{x} \iff a(\bar{x}; u) \to 1$ and $o(\bar{x}; u) \to \theta$, for $u \in H$.

For the sake of simplicity, $a(\bar{x}; \cdot)$ and $o(\bar{x}; \cdot)$ are abbreviated as $a(\cdot)$ and $o(\cdot)$, respectively.

### 3.2. Review the strict Fréchet differentiability of the metric projection operator onto balls in Hilbert spaces

Let $\mathbb{B}$ denotes the unit closed ball in a Hilbert space $H$. For any $r > 0$, $r\mathbb{B}$ denotes the closed ball with radius $r$ and centered at $\theta$. Let $\mathbb{S}$ be the unit sphere in $H$. Then, $r\mathbb{S}$ is the sphere in $H$ with radius $r$ and centered at $\theta$. For any $c \in H$ and $r > 0$, let $\mathbb{B}(c, r)$ denote the closed ball in $H$ with radius $r$ and centered at $c$. In this notation, $\mathbb{B}(\theta, 1) = \mathbb{B}$ and $\mathbb{B}(\theta, r) = r\mathbb{B}$. $\mathbb{S}(c, r)$ denotes the sphere in $H$ with center $c$ and with radius $r$. Let $I_H$ denote the identity mapping in $H$. To study the differentiability of $P_{r\mathbb{B}}$, we need the following representations of $P_{r\mathbb{B}}$.

(a) $\quad P_{r\mathbb{B}}(x) = x$, for any $x \in r\mathbb{B}$;

(b) $\quad P_{r\mathbb{B}}(x) = \frac{r}{\|x\|} x$, for any $x \in H \setminus r\mathbb{B}$.

In [8], the strict Fréchet differentiability of $P_{r\mathbb{B}}$ is obtained in Hilbert spaces. To state those results from [8], we need recall the following notations.

For any $x \in \mathbb{S}(c, r)$, we define two subsets $x^{\uparrow}_{(c,r)}$ and $x^{\downarrow}_{(c,r)}$ in $H \setminus \{\theta\}$ as follows:

(a) $x^{\uparrow}_{(c,r)} = \{v \in H \setminus \{\theta\}:$ there is $\delta > 0$ such that $\|(x + tv) - c\| \geq r$, for all $t \in (0, \delta)\}$;

(b) $x^{\downarrow}_{(c,r)} = \{v \in H \setminus \{\theta\}:$ there is $\delta > 0$ such that $\|(x + tv) - c\| < r$, for all $t \in (0, \delta)\}$;

(c) $x^{\uparrow}_r = x^{\uparrow}_{(\theta,r)}$;

(d) $x^{\downarrow}_r = x^{\downarrow}_{(\theta,r)}$.

**Theorem 3.3 in [8].** *Let $H$ be a Hilbert space. For any $r > 0$, the metric projection $P_{r\mathbb{B}}: H \to r\mathbb{B}$ has the following differentiability properties.*

(i) $\quad P_{r\mathbb{B}}$ is strictly Fréchet differentiable on $r\mathbb{B}^\circ$ satisfying

$$\nabla P_{r\mathbb{B}}(\bar{x}) = I_H, \text{ for every } \bar{x} \in r\mathbb{B}^\circ.$$

That is,

$$\bar{x} \in r\mathbb{B}^\circ \implies \nabla P_{r\mathbb{B}}(\bar{x})(x) = x, \text{ for every } x \in H.$$

(ii) $\quad P_{r\mathbb{B}}$ is strictly Fréchet differentiable on $H \setminus r\mathbb{B}$ such that, for every $\bar{x} \in H \setminus r\mathbb{B}$,

$$\nabla P_{r\mathbb{B}}(\bar{x})(x) = \frac{r}{\|\bar{x}\|} o(\bar{x}; x) = \frac{r}{\|\bar{x}\|}\left(x - \frac{\langle x, \bar{x}\rangle}{\|\bar{x}\|^2}\bar{x}\right), \text{ for every } x \in H.$$

*In particular, we have*

(a) $\quad \nabla P_{r\mathbb{B}}(\bar{x})(x) = \frac{r}{\|\bar{x}\|}x, \text{ if } x \perp \bar{x}, \text{ for } x \in H;$

(b) $\quad \nabla P_{r\mathbb{B}}(\bar{x})(\bar{x}) = \theta.$

(iii) *For the subset* $r\mathbb{S}$, *we have*

(I) $P_{r\mathbb{B}}$ *is Gâteaux directional differentiable on* $r\mathbb{S}$ *satisfying that, for every point* $\bar{x} \in r\mathbb{S}$, *the following representations are satisfied*

(a) $\quad P'_{r\mathbb{B}}(\bar{x})(w) = w - \frac{1}{r^2}\langle \bar{x}, w \rangle x, \quad \text{if } w \in \bar{x}_r^{\uparrow};$

(b) $\quad P'_{r\mathbb{B}}(\bar{x})(\bar{x}) = \theta;$

(c) $\quad P'_{r\mathbb{B}}(\bar{x})(w) = w, \quad \text{if } w \in \bar{x}_r^{\downarrow}.$

(II) $P_{r\mathbb{B}}$ *is not Fréchet differentiable at any* $\bar{x} \in r\mathbb{S}$. *That is,*

$$\nabla P_{r\mathbb{B}}(\bar{x}) \text{ does not exist, for any } \bar{x} \in r\mathbb{S}.$$

### 3.3. Mordukhovich derivatives of the metric projection onto balls in Hilbert spaces

For $r > 0$, we consider the single-valued metric projection operator $P_{r\mathbb{B}}: H \to r\mathbb{B}$. By Lemma 2.1, for any $\bar{x} \in H$, the Mordukhovich derivatives of $P_{r\mathbb{B}}$ at $(\bar{x}, P_{r\mathbb{B}}(\bar{x}))$ is calculated by

$$\widehat{D}^* P_{r\mathbb{B}}(\bar{x})(y) = \left\{ z \in H: \limsup_{(u, P_{r\mathbb{B}}(u)) \xrightarrow{\mathrm{gph} P_{r\mathbb{B}}} (\bar{x}, P_{r\mathbb{B}}(\bar{x}))} \frac{\langle z, u - \bar{x} \rangle - \langle y, P_{r\mathbb{B}}(u) - P_{r\mathbb{B}}(\bar{x}) \rangle}{\|u - \bar{x}\| + \|P_{r\mathbb{B}}(u) - P_{r\mathbb{B}}(\bar{x})\|} \leq 0 \right\}.$$

Since $P_{r\mathbb{B}}: H \to r\mathbb{B}$ is continuous, by Lemma 2.1, for any $\bar{x} \in H$, we have

$$\widehat{D}^* P_{r\mathbb{B}}(\bar{x})(y) = \left\{ z \in H: \limsup_{u \to \bar{x}} \frac{\langle z, u - \bar{x} \rangle - \langle y, P_{r\mathbb{B}}(u) - P_{r\mathbb{B}}(\bar{x}) \rangle}{\|u - \bar{x}\| + \|P_{r\mathbb{B}}(u) - P_{r\mathbb{B}}(\bar{x})\|} \leq 0 \right\}. \quad (3.1)$$

If $P_{r\mathbb{B}}$ is Fréchet differentiable at $\bar{x}$, then, by Theorem 1.38 in [11], the Mordukhovich derivatives of $P_{r\mathbb{B}}$ at $\bar{x}$ is calculated by

$$\widehat{D}^* P_{r\mathbb{B}}(\bar{x})(y) = \{(\nabla P_{r\mathbb{B}}(\bar{x}))(y)\}, \text{ for all } y \in H. \quad (3.2)$$

Now, for $r > 0$, by Theorem 3.3 in [8], and formulas (3.1) and (3.2), we find the Fréchet coderivatives of $P_{r\mathbb{B}}$ at any point in $H$. However, since the calculations in part (III) of the following theorem are considered long and complicated, we only study three cases. Readers who are interested in this topic can similarly study other possible cases.

**Theorem 3.2.** *Let $H$ be a Hilbert space. For any $r > 0$, the Mordukhovich derivatives of $P_{r\mathbb{B}}$ at a point $\bar{x} \in H$ is given below.*

(I) $\quad$ *If* $\bar{x} \in r\mathbb{B}^{\circ}$, *then*

$$\widehat{D}^*P_{r\mathbb{B}}(\bar{x})(y) = y, \text{ for every } y \in H.$$

(II)      If $\bar{x} \in H\setminus r\mathbb{B}$, then

$$\widehat{D}^*P_{r\mathbb{B}}(\bar{x})(y) = \frac{r}{\|\bar{x}\|}\left(y - \frac{\langle y,\bar{x}\rangle}{\|\bar{x}\|^2}\bar{x}\right), \text{ for every } y \in H.$$

(III)      If $\bar{x} \in r\mathbb{S}$, then

    (a)      $\widehat{D}^*P_{r\mathbb{B}}(\bar{x})(\theta) = \{\theta\};$

    (b)      *For any $y \in H\setminus\{\theta\}$, we have*

$$\theta \in \widehat{D}^*P_{r\mathbb{B}}(\bar{x})(y) \iff y = \frac{\langle y,\bar{x}\rangle}{\|\bar{x}\|^2}\bar{x} \text{ with } \langle y,\bar{x}\rangle \leq 0;$$

    (c)      $\widehat{D}^*P_{r\mathbb{B}}(\bar{x})(\bar{x}) = \emptyset.$

*Proof.* Parts (I, II) follow from Theorem 3.3 and (3.2), which is Theorem 1.38 in [11]. So, we only prove part (III).

Proof of (a) in (III). Since $\|\bar{x}\| = r$, we have $P_{r\mathbb{B}}(\bar{x}) = \bar{x}$. It is clear to see

$$\theta \in \widehat{D}^*P_{r\mathbb{B}}(\bar{x})(\theta). \tag{3.3}$$

For $z \in H\setminus\{\theta\}$, by (3.1), we have

$$z \in \widehat{D}^*P_{r\mathbb{B}}(\bar{x})(\theta) \iff \limsup_{u \to \bar{x}} \frac{\langle(z,-\theta),\ (u,P_{r\mathbb{B}}(u))-(\bar{x},P_{r\mathbb{B}}(\bar{x}))\rangle}{\|u-\bar{x}\|+\|P_{r\mathbb{B}}(u)-P_{r\mathbb{B}}(\bar{x})\|} \leq 0$$

$$\iff \limsup_{u \to \bar{x}} \frac{\langle(z,-\theta),\ (u,P_{r\mathbb{B}}(u))-(\bar{x},\bar{x})\rangle}{\|u-\bar{x}\|+\|P_{r\mathbb{B}}(u)-\bar{x}\|} \leq 0$$

$$\iff \limsup_{u \to \bar{x}} \frac{\langle z,\ u-\bar{x}\rangle}{\|u-\bar{x}\|+\|P_{r\mathbb{B}}(u)-\bar{x}\|} \leq 0. \tag{3.4}$$

We take a directional line segment in the limit (3.4) as, $u = (1+\delta)\bar{x}$, for $\delta \downarrow 0$. We have

$$z \in \widehat{D}^*P_{r\mathbb{B}}(\bar{x})(\theta) \iff 0 \geq \limsup_{u \to \bar{x}} \frac{\langle z,\ u-\bar{x}\rangle}{\|u-\bar{x}\|+\|P_{r\mathbb{B}}(u)-\bar{x}\|}$$

$$\geq \limsup_{u \to \bar{x}} \frac{\langle a(z)\bar{x}+o(z),\ (1+\delta)\bar{x}-\bar{x}\rangle}{\|(1+\delta)\bar{x}-\bar{x}\|+\|P_{r\mathbb{B}}((1+\delta)\bar{x})-\bar{x}\|}$$

$$= \limsup_{u \to \bar{x}} \frac{\delta a(z)\|\bar{x}\|^2}{\delta\|\bar{x}\|+\|\bar{x}-\bar{x}\|}$$

$$= \limsup_{u \to \bar{x}} \frac{\delta a(z)\|\bar{x}\|^2}{\delta\|\bar{x}\|}$$

$$= a(z)r.$$

This implies

$$a(z) \leq 0.$$

Similar to the proof of $a(z) \leq 0$, if we take a directional line segment in the limit in (3.4) as, $u = (1 - \delta)\bar{x}$, for $\delta \downarrow 0$, we will obtain $a(z) \geq 0$. This implies that

$$z \in \widehat{D}^* P_{r\mathbb{B}}(\bar{x})(\theta) \implies a(z) = 0. \qquad (3.5)$$

Then, by assumption that $z \in H|\{\theta\}$ and $a(z) = 0$, we must have

$$z = o(z) \neq \theta. \qquad (3.6)$$

We take a directional line segment in the limit in (3.4) as, $u = \bar{x} + \delta o(z)$, for $\delta \downarrow 0$. We have

$$z \in \widehat{D}^* P_{r\mathbb{B}}(\bar{x})(\theta) \iff 0 \geq \limsup_{u \to \bar{x}} \frac{\langle z, u - \bar{x}\rangle}{\|u - \bar{x}\| + \|P_{r\mathbb{B}}(u) - \bar{x}\|}$$

$$\geq \limsup_{\delta \downarrow 0} \frac{\langle o(z), \bar{x} + \delta o(z) - \bar{x}\rangle}{\|\bar{x} + \delta o(z) - \bar{x}\| + \|P_{r\mathbb{B}}(\bar{x} + \delta o(z)) - \bar{x}\|}$$

$$= \limsup_{\delta \downarrow 0} \frac{\langle o(z), \delta o(z)\rangle}{\delta \|o(z)\| + \left\|\frac{r}{\|\bar{x} + \delta o(z)\|}(\bar{x} + \delta o(z)) - \bar{x}\right\|}$$

$$= \limsup_{\delta \downarrow 0} \frac{\delta \|o(z)\|^2}{\delta \|o(z)\| + \left\|\frac{r}{\|\bar{x} + \delta o(z)\|}(\bar{x} + \delta o(z)) - \bar{x}\right\|}. \qquad (3.7)$$

We calculate $\left\|\frac{r}{\|\bar{x} + \delta o(z)\|}(\bar{x} + \delta o(z)) - \bar{x}\right\|$ in the denominator in limit (3.7).

$$\left\|\frac{r}{\|\bar{x} + \delta o(z)\|}(\bar{x} + \delta o(z)) - \bar{x}\right\|^2$$

$$= \left(\frac{r}{\|\bar{x} + \delta o(z)\|} - 1\right)^2 \|\bar{x}\|^2 + \left(\frac{r}{\|\bar{x} + \delta o(z)\|}\right)^2 \|\delta o(z)\|^2$$

$$= \left(\frac{r^2 - \|\bar{x} + \delta o(z)\|^2}{\|\bar{x} + \delta o(z)\|(r + \|\bar{x} + \delta o(z)\|)}\right)^2 + \delta^2 \|o(z)\|^2 \left(\frac{r}{\|\bar{x} + \delta o(z)\|}\right)^2$$

$$= \left(\frac{r^2 - (\|\bar{x}\|^2 + \|\delta o(z)\|^2)}{\|\bar{x} + \delta o(z)\|(r + \|\bar{x} + \delta o(z)\|)}\right)^2 + \delta^2 \|o(z)\|^2 \left(\frac{r}{\|\bar{x} + \delta o(z)\|}\right)^2$$

$$= \left(\frac{\delta^2 \|o(z)\|^2}{\|\bar{x} + \delta o(z)\|(r + \|\bar{x} + \delta o(z)\|)}\right)^2 + \delta^2 \|o(z)\|^2 \left(\frac{r}{\|\bar{x} + \delta o(z)\|}\right)^2$$

$$= \delta^2 \|o(z)\|^2 \left(\left(\frac{\delta \|o(z)\|}{\|\bar{x} + \delta o(z)\|(r + \|\bar{x} + \delta o(z)\|)}\right)^2 + \left(\frac{r}{\|\bar{x} + \delta o(z)\|}\right)^2\right).$$

Hence,

$$\left\|\frac{r}{\|\bar{x} + \delta o(z)\|}(\bar{x} + \delta o(z)) - \bar{x}\right\| = \delta \|o(z)\| \sqrt{\left(\frac{\delta \|o(z)\|}{\|\bar{x} + \delta o(z)\|(r + \|\bar{x} + \delta o(z)\|)}\right)^2 + \left(\frac{r}{\|\bar{x} + \delta o(z)\|}\right)^2}. \qquad (3.8)$$

Substituting (3.8) into (3.7), we have

$$z \in \widehat{D}^* P_{r\mathbb{B}}(\bar{x})(\theta) \iff 0 \geq \limsup_{u \to \bar{x}} \frac{\langle z, u-\bar{x} \rangle}{\|u-\bar{x}\| + \|P_{r\mathbb{B}}(u)-\bar{x}\|}$$

$$\geq \limsup_{\delta \downarrow 0} \frac{\delta \|o(z)\|^2}{\delta \|o(z)\| + \delta \|o(z)\| \sqrt{\left(\frac{\delta \|o(z)\|}{\|\bar{x}+\delta o(z)\|(r+\|\bar{x}+\delta o(z)\|)}\right)^2 + \left(\frac{r}{\|\bar{x}+\delta o(z)\|}\right)^2}}$$

$$= \limsup_{\delta \downarrow 0} \frac{\|o(z)\|^2}{\|o(z)\| + \|o(z)\| \sqrt{\left(\frac{\delta \|o(z)\|}{\|\bar{x}+\delta o(z)\|(r+\|\bar{x}+\delta o(z)\|)}\right)^2 + \left(\frac{r}{\|\bar{x}+\delta o(z)\|}\right)^2}}$$

$$= \frac{\|o(z)\|^2}{\|o(z)\| + \|o(z)\|}$$

$$= \frac{\|o(z)\|}{2}.$$

This implies $\|o(z)\| = 0$, which contradicts to the assumption (3.6). It follows that

$$z \in \widehat{D}^* P_{r\mathbb{B}}(\bar{x})(\theta) \implies o(z) = 0.$$

Hence, by (3.5), this implies

$$z \in \widehat{D}^* P_{r\mathbb{B}}(\bar{x})(\theta) \implies z = \theta.$$

By (3.3), this proves (a) in (III).

Next, we prove the part "$\implies$" of (b) in (III). Let $y \in H$ with $y \neq \theta$. By definition,

$$\theta \in \widehat{D}^* P_{r\mathbb{B}}(\bar{x})(y) \iff \limsup_{u \to \bar{x}} \frac{\langle (\theta, -y), (u, P_{r\mathbb{B}}(u) - (\bar{x}, P_{r\mathbb{B}}(\bar{x})) \rangle}{\|u - \bar{x}\| + \|P_{r\mathbb{B}}(u) - P_{r\mathbb{B}}(\bar{x})\|} \leq 0$$

$$\iff \limsup_{u \to \bar{x}} \frac{-\langle y, P_{r\mathbb{B}}(u) - \bar{x} \rangle}{\|u - \bar{x}\| + \|P_{r\mathbb{B}}(u) - \bar{x}\|} \leq 0. \tag{3.9}$$

We take a directional line segment in the limit in (3.9), $u = (1-\delta)\bar{x}$, for $\delta \downarrow 0$ with $\delta < 1$. Then

$$0 \geq \limsup_{u \to \bar{x}} \frac{-\langle y, P_{r\mathbb{B}}(u) - \bar{x} \rangle}{\|u - \bar{x}\| + \|P_{r\mathbb{B}}(u) - \bar{x}\|}$$

$$\geq \limsup_{\delta \downarrow 0} \frac{-\langle y, P_{r\mathbb{B}}((1-\delta)\bar{x}) - \bar{x} \rangle}{\|(1-\delta)\bar{x} - \bar{x}\| + \|P_{r\mathbb{B}}((1-\delta)\bar{x}) - \bar{x}\|}$$

$$= \limsup_{\delta \downarrow 0} \frac{-\langle a(y)\bar{x} + o(y), -\delta \bar{x} \rangle}{\|(1-\delta)\bar{x} - \bar{x}\| + \|(1-\delta)\bar{x} - \bar{x}\|}$$

$$= \limsup_{\delta \downarrow 0} \frac{\delta a(y) \|\bar{x}\|^2}{2\delta \|\bar{x}\|}$$

$$= \frac{a(y)r}{2}.$$

This implies that
$$\theta \in \widehat{D}^* P_{r\mathbb{B}}(\bar{x})(y) \implies a(y) \leq 0. \tag{3.10}$$

Assume, by the way of contradiction, that
$$o(y) \neq \theta. \tag{3.11}$$

Under the hypothesis (3.11), we take a directional line segment in the limit in (3.9), $u = \bar{x} - \delta o(y)$, for $\delta \downarrow 0$ with $0 < \delta < 1$. We suppose that $\|\bar{x} - \delta o(y)\| > r$. In this case, we have
$$P_{r\mathbb{B}}(\bar{x} - \delta o(y)) = \frac{r}{\|\bar{x} - \delta o(y)\|}(\bar{x} - \delta o(y)).$$

In case, if $\|\bar{x} - \delta o(y)\| \leq r$, then $P_{r\mathbb{B}}(\bar{x} - \delta o(y)) = \bar{x} - \delta o(y)$, then the following calculation will be simpler. So, by (3.11), it follows that

$$0 \geq \limsup_{u \to \bar{x}} \frac{-\langle y, P_{r\mathbb{B}}(u) - \bar{x}\rangle}{\|u - \bar{x}\| + \|P_{r\mathbb{B}}(u) - \bar{x}\|}$$

$$\geq \limsup_{\delta \downarrow 0} \frac{-\langle a(y)\bar{x} + o(y), P_{r\mathbb{B}}(\bar{x} - \delta o(y)) - \bar{x}\rangle}{\|\bar{x} - \delta o(y) - \bar{x}\| + \|P_{r\mathbb{B}}(\bar{x} - \delta o(y)) - \bar{x}\|}$$

$$= \limsup_{\delta \downarrow 0} \frac{-\langle a(y)\bar{x} + o(y), \frac{r}{\|\bar{x} - \delta o(y)\|}(\bar{x} - \delta o(y)) - \bar{x}\rangle}{\delta \|o(y)\| + \left\|\frac{r}{\|\bar{x} - \delta o(y)\|}(\bar{x} - \delta o(y)) - \bar{x}\right\|}$$

$$= \limsup_{\delta \downarrow 0} \frac{-a(y)\left(\frac{r}{\|\bar{x} - \delta o(y)\|} - 1\right)\|\bar{x}\|^2 + \delta \|o(y)\|^2 \frac{r}{\|\bar{x} - \delta o(y)\|}}{\delta \|o(y)\| + \left\|\frac{r}{\|\bar{x} - \delta o(y)\|}(\bar{x} - \delta o(y)) - \bar{x}\right\|}$$

$$= \limsup_{\delta \downarrow 0} \frac{-a(y)r^2\left(\frac{r - \|\bar{x} - \delta o(y)\|}{\|\bar{x} - \delta o(y)\|}\right) + \delta \|o(y)\|^2 \frac{r}{\|\bar{x} - \delta o(y)\|}}{\delta \|o(y)\| + \left\|\frac{r}{\|\bar{x} - \delta o(y)\|}(\bar{x} - \delta o(y)) - \bar{x}\right\|}. \tag{3.12}$$

Similar to (3.8), we calculate $\left\|\frac{r}{\|\bar{x} - \delta o(y)\|}(\bar{x} - \delta o(y)) - \bar{x}\right\|$ to get

$$\left\|\frac{r}{\|\bar{x} - \delta o(y)\|}(\bar{x} - \delta o(y)) - \bar{x}\right\| = \delta \|o(y)\| \sqrt{\left(\frac{\delta \|o(y)\|}{\|\bar{x} - \delta o(y)\|(r + \|\bar{x} - \delta o(y)\|)}\right)^2 + \left(\frac{r}{\|\bar{x} - \delta o(y)\|}\right)^2}. \tag{3.13}$$

Substituting (3.13) into (3.12), we have

$$0 \geq \limsup_{\delta \downarrow 0} \left(\frac{-a(y)r^2\left(\frac{r - \|\bar{x} - \delta o(y)\|}{\|\bar{x} - \delta o(y)\|}\right) + \delta \|o(y)\|^2 \frac{r}{\|\bar{x} - \delta o(y)\|}}{\delta \|o(y)\| + \left\|\frac{r}{\|\bar{x} - \delta o(y)\|}(\bar{x} - \delta o(y)) - \bar{x}\right\|}\right)$$

$$= \limsup_{\delta \downarrow 0} \frac{-a(y)r^2\left(\frac{r^2 - (r^2 + \delta^2 \|o(y)\|^2)}{\|\bar{x} - \delta o(y)\|(r + \|\bar{x} - \delta o(y)\|)}\right) + \delta \|o(y)\|^2 \frac{r}{\|\bar{x} - \delta o(y)\|}}{\delta \|o(y)\| + \delta \|o(y)\| \sqrt{\left(\frac{\delta \|o(y)\|}{\|\bar{x} - \delta o(y)\|(r + \|\bar{x} - \delta o(y)\|)}\right)^2 + \left(\frac{r}{\|\bar{x} - \delta o(y)\|}\right)^2}}$$

$$= \operatorname*{limsup}_{\delta \downarrow 0} \frac{-a(y)r^2 \frac{-\delta^2\|o(y)\|^2}{\|\bar{x}-\delta o(y)\|(r+\|\bar{x}-\delta o(y)\|)} + \delta\|o(y)\|^2 \frac{r}{\|\bar{x}-\delta o(y)\|}}{\delta\|o(y)\|+\delta\|o(y)\| \sqrt{\left(\frac{\delta\|o(y)\|}{\|\bar{x}-\delta o(y)\|(r+\|\bar{x}-\delta o(y)\|)}\right)^2 + \left(\frac{r}{\|\bar{x}-\delta o(y)\|}\right)^2}}$$

$$= \operatorname*{limsup}_{\delta \downarrow 0} \frac{-a(y)r^2 \frac{-\delta\|o(y)\|^2}{\|\bar{x}-\delta o(y)\|(r+\|\bar{x}-\delta o(y)\|)} + \|o(y)\|^2 \frac{r}{\|\bar{x}-\delta o(y)\|}}{\|o(y)\|+\|o(y)\| \sqrt{\left(\frac{\delta\|o(y)\|}{\|\bar{x}-\delta o(y)\|(r+\|\bar{x}-\delta o(y)\|)}\right)^2 + \left(\frac{r}{\|\bar{x}-\delta o(y)\|}\right)^2}}$$

$$= \frac{0 + \|o(y)\|^2 \frac{r}{\|\bar{x}\|}}{\|o(y)\|+\|o(y)\|}$$

$$= \frac{\|o(y)\|}{2}.$$

This implies $o(y) = \theta$, which contradicts to the assumption (3.11). Hence, we obtain that

$$\theta \in \widehat{D}^* P_{r\mathbb{B}}(\bar{x})(y) \implies o(y) = \theta. \tag{3.14}$$

By (3.10) and (3.14) together, we have

$$\theta \in \widehat{D}^* P_{r\mathbb{B}}(\bar{x})(y) \implies y = a(y)\bar{x} \text{ with for } a(y) \le 0. \tag{3.15}$$

Then, we prove the part "$\Longleftarrow$" of (b) in (III). By part (a), it is proved that (3.15) holds for $a(y) = 0$. So, to prove the part "$\Longleftarrow$" of (b), we only need to show that

$$y = a(y)\bar{x} \text{ with } a(y) < 0 \implies \theta \in \widehat{D}^* P_{r\mathbb{B}}(\bar{x})(y). \tag{3.16}$$

Let $y = a(y)\bar{x}$, with $a(y) < 0$. We calculate

$$\operatorname*{limsup}_{u \to \bar{x}} \frac{\langle (\theta, -a(y)\bar{x}),\ (u, P_{r\mathbb{B}}(u) - (\bar{x}, P_{r\mathbb{B}}(\bar{x})) \rangle}{\|u-\bar{x}\| + \| P_{r\mathbb{B}}(u) - P_{r\mathbb{B}}(\bar{x})\|}$$

$$= \operatorname*{limsup}_{u \to \bar{x}} \frac{\langle (\theta, -a(y)\bar{x}),\ (u, P_{r\mathbb{B}}(u) - (\bar{x}, \bar{x}) \rangle}{\|u-\bar{x}\| + \| P_{r\mathbb{B}}(u) - \bar{x}\|}$$

$$= \operatorname*{limsup}_{u \to \bar{x}} \frac{\langle -a(y)\bar{x},\ P_{r\mathbb{B}}(u) - \bar{x} \rangle}{\|u-\bar{x}\| + \| P_{r\mathbb{B}}(u) - \bar{x}\|}$$

$$= \operatorname*{limsup}_{u \to \bar{x}} \frac{a(y)(\|\bar{x}\|^2 - \langle \bar{x},\ P_{r\mathbb{B}}(u) \rangle)}{\|u-\bar{x}\| + \| P_{r\mathbb{B}}(u) - \bar{x}\|}$$

$$\le \operatorname*{limsup}_{u \to \bar{x}} \frac{a(y)(\|\bar{x}\|^2 - \|\bar{x}\|\|P_{r\mathbb{B}}(u)\|)}{\|u-\bar{x}\| + \| P_{r\mathbb{B}}(u) - \bar{x}\|}$$

$$= \operatorname*{limsup}_{u \to \bar{x}} \frac{a(y)r(r - \|P_{r\mathbb{B}}(u)\|)}{\|u-\bar{x}\| + \| P_{r\mathbb{B}}(u) - \bar{x}\|}$$

$$\le \operatorname*{limsup}_{u \to \bar{x}} \frac{0}{\|u-\bar{x}\| + \| P_{r\mathbb{B}}(u) - \bar{x}\|}$$

$$= 0.$$

By definition, this implies $\theta \in \widehat{D}^* P_{r\mathbb{B}}(\bar{x})(y)$, for $y = a(y)\bar{x}$ with $a(y) < 0$, which proves (3.16). By the orthogonal decomposition of $y$ with respect to $\bar{x}$, $a(y) = \frac{\langle y, \bar{x}\rangle}{\|\bar{x}\|^2}$. It follows that $a(y) < 0$ if and only if $\langle y, \bar{x}\rangle < 0$. Then, by (3.15) and (3.16), part (b) is proved.

Proof of (c). Take an arbitrary $z \in H$, by definition, we have

$$z \in \widehat{D}^* P_{r\mathbb{B}}(\bar{x})(\bar{x}) \iff \limsup_{u \to \bar{x}} \frac{\langle (z,-\bar{x}), (u, P_{r\mathbb{B}}(u) - (\bar{x}, P_{r\mathbb{B}}(\bar{x})))\rangle}{\|u-\bar{x}\| + \|P_{r\mathbb{B}}(u) - P_{r\mathbb{B}}(\bar{x})\|} \leq 0$$

$$\iff \limsup_{u \to \bar{x}} \frac{\langle (z,-\bar{x}), (u, P_{r\mathbb{B}}(u) - (\bar{x}, \bar{x}))\rangle}{\|u-\bar{x}\| + \|P_{r\mathbb{B}}(u) - \bar{x}\|} \leq 0$$

$$\iff \limsup_{u \to \bar{x}} \frac{\langle z, u-\bar{x}\rangle - \langle \bar{x}, P_{r\mathbb{B}}(u) - \bar{x}\rangle}{\|u-\bar{x}\| + \|P_{r\mathbb{B}}(u) - \bar{x}\|} \leq 0. \tag{3.17}$$

We take a directional line segment in the limit in (3.17), $u = (1-\delta)\bar{x}$, for $\delta \downarrow 0$ with $\delta < 1$. Then

$$0 \geq \limsup_{u \to \bar{x}} \frac{\langle z, u-\bar{x}\rangle - \langle \bar{x}, P_{r\mathbb{B}}(u) - \bar{x}\rangle}{\|u-\bar{x}\| + \|P_{r\mathbb{B}}(u) - \bar{x}\|}$$

$$\geq \limsup_{\delta \downarrow 0} \frac{\langle z, (1-\delta)\bar{x}-\bar{x}\rangle - \langle \bar{x}, (1-\delta)\bar{x}-\bar{x}\rangle}{\|(1-\delta)\bar{x}-\bar{x}\| + \|(1-\delta)\bar{x}-\bar{x}\|}$$

$$= \limsup_{\delta \downarrow 0} \frac{\langle a(z)\bar{x} + o(z), -\delta\bar{x}\rangle - \langle \bar{x}, -\delta\bar{x}\rangle}{2\delta\|\bar{x}\|}$$

$$= \frac{-a(z)\|\bar{x}\|^2 + \|\bar{x}\|^2}{2\|\bar{x}\|}$$

$$= \frac{(1-a(z))\|\bar{x}\|}{2}.$$

It implies that

$$z \in \widehat{D}^* P_{r\mathbb{B}}(\bar{x})(\bar{x}) \implies a(z) \geq 1. \tag{3.18}$$

If we take a directional line segment in the limit in (3.17) as, $u = (1+\delta)\bar{x}$, for $\delta \downarrow 0$. Then

$$0 \geq \limsup_{u \to \bar{x}} \frac{\langle z, u-\bar{x}\rangle - \langle \bar{x}, P_{r\mathbb{B}}(u) - \bar{x}\rangle}{\|u-\bar{x}\| + \|P_{r\mathbb{B}}(u) - \bar{x}\|}$$

$$\geq \limsup_{\delta \downarrow 0} \frac{\langle z, (1+\delta)\bar{x}-\bar{x}\rangle - \langle \bar{x}, P_{r\mathbb{B}}((1+\delta)\bar{x}) - \bar{x}\rangle}{\|(1+\delta)\bar{x}-\bar{x}\| + \|P_{r\mathbb{B}}((1+\delta)\bar{x}) - \bar{x}\|}$$

$$= \limsup_{\delta \downarrow 0} \frac{\langle a(z)\bar{x} + o(z), \delta\bar{x}\rangle - \langle \bar{x}, \bar{x}-\bar{x}\rangle}{\delta\|\bar{x}\| + \|\bar{x}-\bar{x}\|}$$

$$= \frac{a(z)\|\bar{x}\|^2 - 0}{\|\bar{x}\|}$$

$$= a(z)r.$$

This implies that

$$z \in \widehat{D}^* P_{rB}(\bar{x})(\bar{x}) \implies a(z) \leq 0. \tag{3.19}$$

It is clear that there is no point $z \in H$, which simultaneously satisfies both (3.18) and (3.19). This implies that $\widehat{D}^* P_{rB}(\bar{x})(\bar{x}) = \emptyset$, which proves (c). □

## 4. Mordukhovich derivatives of the metric projection onto the positive cone in $\mathbb{R}^n$

### 4.1 Review strict Fréchet differentiability of the metric projection onto the positive cone in $\mathbb{R}^n$

Let $\mathbb{R}^n$ denote the $n$-dimensional Euclidean space with origin $\theta$. We recall some notations and concepts in $\mathbb{R}^n$ from [8], which will be used in this section. We define

$$\Delta \mathbb{R}^n = \{x = (x_1, x_2, \ldots, x_n) \in \mathbb{R}^n : x_k = 0, \text{ for at least one } k \in \{1, 2, \ldots, n\}\}.$$

Let $K$ denote the positive cone of $\mathbb{R}^n$ defined by

$$K = \{x = (x_1, x_2, \ldots, x_n) \in \mathbb{R}^n : x_i \geq 0, i = 1, 2, \ldots, n\}.$$

$K$ is a pointed closed and convex cone in $\mathbb{R}^n$. The interior of $K$ is denoted by $K^o$, which is a nonempty subset of $K$ satisfying

$$K^o = \{x = (x_1, x_2, \ldots, x_n) \in \mathbb{R}^n : x_i > 0, i = 1, 2, \ldots, n\}.$$

The boundary of $K$ is denoted by $\partial K$ such that

$$\partial K = \{x = (x_1, x_2, \ldots, x_n) \in K : \prod_{i=1}^n x_i = 0\}.$$

Then, we see that $\partial K \subseteq \Delta \mathbb{R}^n$. The negative cone of $\mathbb{R}^n$ is $-K$ satisfying $K \cap (-K) = \{\theta\}$. $-K$ is also a pointed closed and convex cone in $\mathbb{R}^n$. We can similarly define the interior $(-K)^o$ and the boundary $\partial(-K)$ of $-K$. For any $x = (x_1, x_2, \ldots, x_n) \in \mathbb{R}^n$, we define three subsets of the set $\{1, 2, \ldots, n\}$ with respect to the given $x$ by

$$x^+ = \{i \in \{1, 2, \ldots, n\} : x_i > 0\},$$

$$x^- = \{i \in \{1, 2, \ldots, n\} : x_i < 0\},$$

and

$$\dot{x} = \{i \in \{1, 2, \ldots, n\} : x_i = 0\}.$$

Define

$$\widehat{K} = \{x = (x_1, x_2, \ldots, x_n) \in \mathbb{R}^n : x^+ \neq \emptyset, x^- \neq \emptyset \text{ and } \dot{x} = \emptyset\}.$$

Then, for any $x = (x_1, x_2, \ldots, x_n) \in \mathbb{R}^n$, we have

(a) $x^+ \cup x^- \cup \dot{x} = \{1, 2, \ldots, n\}$;
(b) $x \in K \iff x^- = \emptyset$;

(c) $x \in K^o \iff x^- = \emptyset$ and $\dot{x} = \emptyset$;
(d) $x \in \partial K \iff x^- = \emptyset$ and $\dot{x} \neq \emptyset$;
(e) $x \in \widehat{K} \iff x^+ \neq \emptyset, x^- \neq \emptyset$ and $\dot{x} = \emptyset$;
(f) $x \in \Delta \mathbb{R}^n \iff \dot{x} \neq \emptyset$.

**Lemma 4.1 in [8].** *Let $K$ be the positive cone of $\mathbb{R}^n$ with negative cone $-K$. $P_K$ has the following properties.*

(a) *For any $x \in \mathbb{R}^n$, $P_K(x)$ is represented as follows*

$$P_K(x)_i = \begin{cases} x_i, & \text{if } i \in x^+, \\ 0, & \text{if } i \notin x^+, \end{cases} \text{ for } i = 1, 2, \ldots, n. \tag{4.1}$$

*In particular, we have*

(i) $P_K(x) = x$, for any $x \in K$;
(ii) $P_K(x) = \theta$, for any $x \in -K$;
(iii) $P_K(x) \in \partial K$, for any $x \in \mathbb{R}^n \setminus K$;
(iv) $P_K(x) \in \partial K \setminus \{\theta\}$, for any $x \in \widehat{K}$.

(b) $P_K$ *is positive homogeneous. For any $x \in \mathbb{R}^n$,*

$$P_K(\lambda x) = \lambda P_K(x), \text{ for any } \lambda \geq 0.$$

For any given fixed $x \in \widehat{K}$, we define a mapping $b(x; \cdot) \colon \mathbb{R}^n \to \mathbb{R}^n$, for any $w \in \mathbb{R}^n$, by

$$(b(x; w))_i = \begin{cases} w_i, & \text{if } i \in x^+, \\ 0, & \text{if } i \notin x^+, \end{cases} \text{ for } i = 1, 2, \ldots, n. \tag{4.2}$$

**Lemma 4.2 in [8].** *For any fixed $x \in \widehat{K}$, as defined in (4.2), $b(x; \cdot)$ satisfies*

(a) $b(x; \cdot) \colon \mathbb{R}^n \to \mathbb{R}^n$ *is a linear and continuous mapping;*
(b) $b(x; w) \in \Delta \mathbb{R}^n$, *for any $w \in \mathbb{R}^n$.*

For any given $x \in \Delta \mathbb{R}^n$, we define a mapping $d(x; \cdot) \colon \mathbb{R}^n \to \mathbb{R}^n$, for any $w \in \mathbb{R}^n$, by

$$\begin{aligned}(d(x; w))_i &= w_i, & \text{for } i \in x^+, \\ (d(x; w))_i &= 0, & \text{for } i \in x^-, \end{aligned}$$

and $\quad (d(x; w))_i = \begin{cases} w_i, & \text{if } w_i > 0, \\ 0, & \text{if } w_i \leq 0, \end{cases} \text{ for } i \in \dot{x}. \tag{4.3}$

**Lemma 4.3 in [8].** *For any fixed $x \in \Delta \mathbb{R}^n$, as defined in (4.3), $d(x; \cdot)$ satisfies*

(a) $d(x; \cdot)$ *is a non-liner mapping from $\mathbb{R}^n$ to $\mathbb{R}^n$;*
(b) $d(\theta; w) = P_K(w)$, *for any $w \in \mathbb{R}^n$.*

**Theorem 4.4 in [8].** *Let $K$ be the positive cone of $\mathbb{R}^n$ with negative cone $-K$. Then, the projection operator $P_K$ has the following Fréchet differentiability.*

(i)      $P_K$ is strict Fréchet differentiable on $K^o$ satisfying $\nabla P_K(x) = I_{\mathbb{R}^n}$, for any $x \in K^o$,

$$\nabla P_K(x)(y) = y, \text{ for any } y \in \mathbb{R}^n;$$

(ii)      $P_K$ is strict Fréchet differentiable on $(-K)^o$ satisfying $\nabla P_K(x) = \theta$, for any $x \in (-K)^o$,

$$\nabla P_K(x)(y) = \theta, \text{ for any } y \in \mathbb{R}^n;$$

(iii)      $P_K$ is strictly Fréchet differentiable on $\widehat{K}$ such that, for any $x \in \widehat{K}$,

$$\nabla P_K(x)(y) = b(x; y), \text{ for any } y \in \mathbb{R}^n;$$

(iv)      For subset $\Delta \mathbb{R}^n$, we have

        (a)    $P_K$ is Gâteaux directionally differentiable on $\Delta \mathbb{R}^n$ such that, for any $x \in \Delta \mathbb{R}^n$,

$$P'_K(x)(w) = d(x; w), \text{ for any } w \in \mathbb{R}^n \setminus \{\theta\}.$$

        (b)    $P_K$ is not Fréchet differentiable at any point in $\Delta \mathbb{R}^n$, that is,

$$\nabla P_K(x) \text{ does not exist, for any } x \in \Delta \mathbb{R}^n.$$

### 4.2 Mordukhovich derivatives of the metric projection onto the positive cone in $\mathbb{R}^n$

Since $P_K : \mathbb{R}^n \to K$ is also continuous, by Lemma 2.1, for any $\bar{x} \in H$, the Mordukhovich derivatives of $P_K$ at point $(\bar{x}, P_K(\bar{x}))$ is calculated by

$$\widehat{D}^* P_K(\bar{x})(y) = \left\{ z \in H : \limsup_{u \to \bar{x}} \frac{\langle z, u-\bar{x} \rangle - \langle y, P_K(u) - P_K(\bar{x}) \rangle}{\|u-\bar{x}\| + \|P_K(u) - P_K(\bar{x})\|} \leq 0 \right\}. \quad (4.4)$$

In particular, if $P_K$ is Fréchet differentiable at $\bar{x}$, then, by Theorem 1.38 in [11], the Mordukhovich derivatives of $P_K$ at $\bar{x}$ is calculated by

$$\widehat{D}^* P_K(\bar{x})(y) = \{(\nabla P_K(\bar{x}))(y)\}, \text{ for all } y \in H. \quad (4.5)$$

By applying Theorem 4.4 in [8] and (4.5), we find the Mordukhovich derivatives of $P_K$.

**Theorem 4.1**. *Let $K$ be the positive cone of $\mathbb{R}^n$ with negative cone $-K$. Then, the Mordukhovich derivatives of the projection operator $P_K$ have the following representations.*

(i)      *For any $\bar{x} \in K^o$, we have*

$$\widehat{D}^* P_K(\bar{x})(y) = y, \text{ for any } y \in \mathbb{R}^n;$$

(ii)      *For any $\bar{x} \in (-K)^o$, we have*

$$\widehat{D}^* P_K(\bar{x})(y) = \theta, \text{ for any } y \in \mathbb{R}^n;$$

(iii) *For any $\bar{x} \in \widehat{K}$, we have*

$$\widehat{D}^* P_K(\bar{x})(y) = b(\bar{x}; y), \text{ for any } y \in \mathbb{R}^n;$$

(iv) *For any $\bar{x} \in \Delta \mathbb{R}^n$, we have*

(a) $\widehat{D}^* P_K(\bar{x})(\theta) = \{\theta\}$;
(b) *For any $y \in \mathbb{R}^n$,*

$$y^- \cap \dot{\bar{x}} \neq \emptyset \implies \lambda y \notin \widehat{D}^* P_K(\bar{x})(y), \text{ for any } \lambda < 1.$$

*In particular,*

$$y^- \cap \dot{\bar{x}} \neq \emptyset \implies \theta \notin \widehat{D}^* P_K(\bar{x})(y);$$

(c) *Let $\bar{x} \in \Delta \mathbb{R}^n \setminus \{\theta\}$. For any $z \in \mathbb{R}^n$, we have,*

$$z \neq \theta \implies z \notin \widehat{D}^* P_K(\bar{x})(\bar{x}), \tag{4.6}$$

$$\bar{x}^+ = \emptyset \iff \widehat{D}^* P_K(\bar{x})(\bar{x}) = \{\theta\}, \tag{4.7}$$

*and*

$$\bar{x}^+ \neq \emptyset \implies \widehat{D}^* P_K(\bar{x})(\bar{x}) = \emptyset. \tag{4.8}$$

*Proof.* Parts (i), (ii) and (iii) follow from Theorem 4.4 in [8] and (4.5) immediately. We only prove (iv). We prove (a) of (iv). For any $\bar{x} \in \Delta \mathbb{R}^n$, it is easy to see that

$$\theta \in \widehat{D}^* P_K(\bar{x})(\theta). \tag{4.9}$$

For any $z \in \mathbb{R}^n$ with $z \neq \theta$, we take a directional line segment $u = \bar{x} + tz$, for $t \downarrow 0$. By the non-expansiveness of the metric projection operator in Hilbert spaces, we have

$$\limsup_{u \to \bar{x}} \frac{\langle z, u - \bar{x} \rangle - \langle \theta, P_K(u) - P_K(\bar{x}) \rangle}{\|u - \bar{x}\| + \|P_K(u) - P_K(\bar{x})\|}$$

$$\geq \limsup_{t \downarrow 0} \frac{\langle z, \bar{x} + tz - \bar{x} \rangle}{\|\bar{x} + tz - \bar{x}\| + \|P_K(\bar{x} + tz) - P_K(\bar{x})\|}$$

$$= \limsup_{t \downarrow 0} \frac{t\|z\|^2}{t\|z\| + \|P_K(\bar{x} + tz) - P_K(\bar{x})\|}$$

$$\geq \limsup_{t \downarrow 0} \frac{t\|z\|^2}{t\|z\| + \|\bar{x} + tz - \bar{x}\|} \quad \text{(Here, non-expansiveness of } P_K \text{ in Hilbert space } H\text{)}$$

$$= \limsup_{t \downarrow 0} \frac{t\|z\|^2}{2t\|z\|}$$

$$= \frac{\|z\|}{2} > 0.$$

This implies

$$z \notin \widehat{D}^* P_K(\bar{x})(\theta), \text{ for any } z \in \mathbb{R}^n \text{ with } z \neq \theta. \tag{4.10}$$

Then, part (a) of (iv) follows from (4.9) and (4.10).

Next, we prove (b) of (iv). Notice that, for any $\bar{x} \in \Delta\mathbb{R}^n$, we have $\dot{\bar{x}} \neq \emptyset$. Let $y \in \mathbb{R}^n$. Suppose $y^- \cap \dot{\bar{x}} \neq \emptyset$. Then, there is $j \in \{1, 2, \ldots, n\}$ such that

$$y_j < 0 \quad \text{and} \quad \bar{x}_j = 0.$$

We take a directional segment $u$ with $u_i = \bar{x}_i$, for $i \neq j$ and $u_j = t$, for $t \downarrow 0$. This implies

$$(u - \bar{x})_i = P_K(u)_i - P_K(\bar{x})_i = \begin{cases} 0, & \text{if } i \neq j, \\ t, & \text{if } i = j, \end{cases} \text{ for } i = 1, 2, \ldots, n.$$

Let $\lambda < 1$. We have

$$\limsup_{u \to \bar{x}} \frac{\langle \lambda y, u - \bar{x} \rangle - \langle y, P_K(u) - P_K(\bar{x}) \rangle}{\|u - \bar{x}\| + \|P_K(u) - P_K(\bar{x})\|}$$

$$\geq \limsup_{t \downarrow 0} \frac{\langle \lambda y, u - \bar{x} \rangle - \langle y, P_K(u) - P_K(\bar{x}) \rangle}{\|u - \bar{x}\| + \|P_K(u) - P_K(\bar{x})\|}$$

$$= \limsup_{t \downarrow 0} \frac{t\lambda y_j - t y_j}{2t}$$

$$= \limsup_{t \downarrow 0} \frac{t(1-\lambda)|y_j|}{2t}$$

$$= \frac{(1-\lambda)|y_j|}{2} > 0, \text{ for any } \lambda < 1.$$

This implies that, for any $y \in \mathbb{R}^n$ with $y^- \cap \dot{\bar{x}} \neq \emptyset$, we have $\lambda y \notin \widehat{D}^* P_K(\bar{x})(y)$, for any $\lambda < 1$. This proves (b) of (iv).

Proof of (c) of (iv). At first, we prove (4.6) in part (c) of (iv). For any $z \in \mathbb{R}^n$ with $z \neq \theta$. The proof of (4.6) is divided to the following two cases.

Case 1. Suppose $(z^+ \cup z^-) \cap \dot{\bar{x}} \neq \emptyset$. Then there is $j \in \{1, 2, \ldots, n\}$ such that

$$\bar{x}_j = 0 \quad \text{and} \quad z_j \neq 0$$

We take a directional segment $u$ with $u_i = \bar{x}_i$, for $i \neq j$ and $u_j = tz_j$, for $t \downarrow 0$. This implies

$$u_i - \bar{x}_i = \begin{cases} 0, & \text{if } i \neq j, \\ tz_j, & \text{if } i = j, \end{cases} \text{ for } i = 1, 2, \ldots, n,$$

and 
$$P_K(u)_i - P_K(\bar{x})_i = \begin{cases} 0, & \text{if } i \neq j, \\ (P_K(u))_j, & \text{if } i = j, \end{cases} \text{ for } i = 1, 2, \ldots, n.$$

This implies that

$$\langle \bar{x}, P_K(u) - P_K(\bar{x}) \rangle$$
$$= \sum_{i \neq j} \bar{x}_i (P_K(u)_i - P_K(\bar{x})_i) + \bar{x}_j (P_K(u)_j - P_K(\bar{x})_j)$$
$$= \sum_{i \neq j} \bar{x}_i 0 + 0 (P_K(u)_j)$$
$$= 0.$$

We have

$$\limsup_{u \to \bar{x}} \frac{\langle z, u - \bar{x} \rangle - \langle \bar{x}, P_K(u) - P_K(\bar{x}) \rangle}{\|u - \bar{x}\| + \|P_K(u) - P_K(\bar{x})\|}$$
$$\geq \limsup_{t \downarrow 0} \frac{\langle z, u - \bar{x} \rangle - 0}{\|u - \bar{x}\| + \|P_K(u) - P_K(\bar{x})\|}$$
$$\geq \limsup_{t \downarrow 0} \frac{\langle z, u - \bar{x} \rangle - 0}{\|u - \bar{x}\|}$$
$$= \limsup_{t \downarrow 0} \frac{t z_j^2}{t |z_j|}$$
$$= |z_j| > 0.$$

This implies that, for any $z \in \mathbb{R}^n$, if $(z^+ \cup z^-) \cap \dot{\bar{x}} \neq \emptyset$, then $z \notin \widehat{D}^* P_K(\bar{x})(\bar{x})$.

Case 2. Suppose $(z^+ \cup z^-) \cap \dot{\bar{x}} = \emptyset$. This implies $\dot{\bar{x}} \subseteq \dot{z}$. Since $z \neq \theta$, there is $j \in \{1, 2, \ldots, n\}$ such that

$$\bar{x}_j \neq 0 \text{ and } z_j \neq 0$$

Subcase 1. Suppose that $\bar{x}_j > 0$ and $z_j \neq 0$. In this case, we take a directional line segment $u$ with $u_i = \bar{x}_i$, for $i \neq j$ and

$$u_j = \begin{cases} \bar{x}_j - t z_j, & \text{if } z_j(\bar{x}_j - z_j) > 0, \\ \bar{x}_j + t z_j, & \text{if } z_j(\bar{x}_j - z_j) < 0, \end{cases} \text{ for } t \downarrow 0.$$

If $0 < t < \frac{\bar{x}_j}{|z_j|}$, then $\bar{x}_j \mp t z_j > 0$. This implies

$$u_i - \bar{x}_i = P_K(u)_i - P_K(\bar{x})_i = \begin{cases} 0, & \text{if } i \neq j, \\ \begin{cases} -t z_j, & \text{if } z_j(\bar{x}_j - z_j) > 0, \\ t z_j, & \text{if } z_j(\bar{x}_j - z_j) < 0, \end{cases} & \text{if } i = j, \end{cases} \text{ for } i = 1, 2, \ldots, n,$$

and
$$\langle \bar{x}, P_K(u) - P_K(\bar{x}) \rangle$$
$$= \sum_{i \neq j} \bar{x}_i (P_K(u)_i - P_K(\bar{x})_i) + \bar{x}_j (P_K(u)_j - P_K(\bar{x})_j)$$
$$= \sum_{i \neq j} \bar{x}_i 0 + \bar{x}_j (P_K(u)_j - P_K(\bar{x})_j)$$

$$= \begin{cases} -t\bar{x}_j z_j, & \text{if } z_j(\bar{x}_j - z_j) > 0, \\ t\bar{x}_j z_j, & \text{if } z_j(\bar{x}_j - z_j) < 0. \end{cases}$$

If $z_j(\bar{x}_j - z_j) > 0$, we have

$$\limsup_{u \to \bar{x}} \frac{\langle z, u-\bar{x} \rangle - \langle \bar{x}, P_K(u) - P_K(\bar{x}) \rangle}{\|u-\bar{x}\| + \|P_K(u) - P_K(\bar{x})\|}$$

$$\geq \limsup_{t \downarrow 0, t < \frac{\bar{x}_j}{|z_j|}} \frac{\langle z, u-\bar{x} \rangle - \langle \bar{x}, P_K(u) - P_K(\bar{x}) \rangle}{\|u-\bar{x}\| + \|P_K(u) - P_K(\bar{x})\|}$$

$$= \limsup_{t \downarrow 0, t < \frac{\bar{x}_j}{|z_j|}} \frac{-tz_j^2 - (-t\bar{x}_j z_j)}{t|z_j| + |\bar{x}_j|}$$

$$= \limsup_{t \downarrow 0} \frac{tz_j(\bar{x}_j - z_j)}{2t|z_j|}$$

$$= \frac{z_j(\bar{x}_j - z_j)}{2|z_j|} > 0.$$

If $z_j(\bar{x}_j - z_j) < 0$, we have

$$\limsup_{u \to \bar{x}} \frac{\langle z, u-\bar{x} \rangle - \langle \bar{x}, P_K(u) - P_K(\bar{x}) \rangle}{\|u-\bar{x}\| + \|P_K(u) - P_K(\bar{x})\|}$$

$$\geq \limsup_{t \downarrow 0, t < \frac{\bar{x}_j}{|z_j|}} \frac{\langle z, u-\bar{x} \rangle - \langle \bar{x}, P_K(u) - P_K(\bar{x}) \rangle}{\|u-\bar{x}\| + \|P_K(u) - P_K(\bar{x})\|}$$

$$= \limsup_{t \downarrow 0, t < \frac{\bar{x}_j}{|z_j|}} \frac{tz_j^2 - (t\bar{x}_j z_j)}{t|z_j| + |\bar{x}_j|}$$

$$= \limsup_{t \downarrow 0} \frac{-tz_j(\bar{x}_j - z_j)}{2t|z_j|}$$

$$= \frac{-z_j(\bar{x}_j - z_j)}{2|z_j|} > 0.$$

This implies that, for any $z \in \mathbb{R}^n$, if $\bar{x}_j > 0$ and $z_j \neq 0$, then $z \notin \widehat{D}^* P_K(\bar{x})(\bar{x})$.

Subcase 2. Suppose that $\bar{x}_j < 0$ and $z_j \neq 0$. In this case, we take a directional segment $u$ with $u_i = \bar{x}_i$, for $i \neq j$ and $u_j = \bar{x}_j + tz_j$, for $t \downarrow 0$. Then, if $0 < t < \frac{|\bar{x}_j|}{|z_j|}$, we have $u_j = \bar{x}_j + tz_j < 0$. This implies, for $0 < t < \frac{|\bar{x}_j|}{|z_j|}$, we get

$$u_i - \bar{x}_i = \begin{cases} 0, & \text{if } i \neq j, \\ tz_j, & \text{if } i = j, \end{cases} \text{ for } i = 1, 2, \ldots, n,$$

and
$$P_K(u) - P_K(\bar{x}) = \theta.$$

We have

$$\limsup_{u \to \bar{x}} \frac{\langle z, u-\bar{x} \rangle - \langle \bar{x}, P_K(u) - P_K(\bar{x}) \rangle}{\|u-\bar{x}\| + \|P_K(u) - P_K(\bar{x})\|}$$

$$\geq \limsup_{t \downarrow 0, t < \frac{|\bar{x}_j|}{|z_j|}} \frac{\langle z, u-\bar{x} \rangle - \langle \bar{x}, P_K(u) - P_K(\bar{x}) \rangle}{\|u-\bar{x}\| + \|P_K(u) - P_K(\bar{x})\|}$$

$$= \limsup_{t \downarrow 0, t < \frac{|\bar{x}_j|}{|z_j|}} \frac{tz_j^2}{t|z_j| + 0}$$

$$= |z_j| > 0.$$

This implies that, for any $z \in \mathbb{R}^n$, if $\bar{x}_j < 0$ and $z_j \neq 0$, then $z \notin \widehat{D}^* P_K(\bar{x})(\bar{x})$.

Then, (4.6) is proved. Next, we prove (4.7) in part (c) of (iv). For any $\bar{x} \in \Delta \mathbb{R}^n \setminus \{\theta\}$ with $\dot{\bar{x}} \neq \emptyset$ and $\bar{x} \neq \theta$. Suppose $\bar{x}^+ = \emptyset$. Since $\bar{x} \neq \theta$, it implies $\bar{x}^- \neq \emptyset$. Let $\delta = \min\{|\bar{x}_i| : \bar{x}_i < 0, \text{ for } i = 1, 2, \ldots, n\}$. By $\bar{x}^- \neq \emptyset$, we have $\delta > 0$. Then, for any $u \in \mathbb{R}^n$,

$$|u| < \frac{\delta}{2} \implies P_K(u)_i - P_K(\bar{x})_i = 0 - 0 = 0, \text{ if } \bar{x}_i < 0, \text{ for } i = 1, 2, \ldots, n.$$

By the assumption that $\bar{x}^+ = \emptyset$. This implies that, if $|u| < \frac{\delta}{2}$, then

$$\langle \bar{x}, P_K(u) - P_K(\bar{x}) \rangle$$

$$= \sum_{i \in \dot{\bar{x}}} \bar{x}_i (P_K(u)_i - P_K(\bar{x})_i) + \sum_{i \in \bar{x}^-} \bar{x}_i (P_K(u)_i - P_K(\bar{x})_i)$$

$$= \sum_{i \in \dot{\bar{x}}} 0 (P_K(u)_i - P_K(\bar{x})_i) + \sum_{i \in \bar{x}^-} \bar{x}_i 0$$

$$= 0.$$

It follows that

$$\limsup_{u \to \bar{x}} \frac{\langle \theta, u-\bar{x} \rangle - \langle \bar{x}, P_K(u) - P_K(\bar{x}) \rangle}{\|u-\bar{x}\| + \|P_K(u) - P_K(\bar{x})\|}$$

$$= \limsup_{u \to \bar{x}, |u| < \frac{\delta}{2}} \frac{\langle \bar{x}, P_K(u) - P_K(\bar{x}) \rangle}{\|u-\bar{x}\| + \|P_K(u) - P_K(\bar{x})\|}$$

$$= 0.$$

This implies that, for any $\bar{x} \in \Delta\mathbb{R}^n$ with $\bar{x}^+ = \emptyset$, we have

$$\theta \in \widehat{D}^* P_K(\bar{x})(\bar{x}).$$

This proves

$$\bar{x}^+ = \emptyset \implies \theta \in \widehat{D}^* P_K(\bar{x})(\bar{x}). \tag{4.11}$$

Next, suppose $\bar{x}^+ \neq \emptyset$. Then, there is $j \in \{1, 2, \ldots, n\}$ such that

$$\bar{x}_j > 0$$

In this case, we take a directional segment $u$ with $u_i = \bar{x}_i$, for $i \neq j$ and $u_j = \bar{x}_j - t$, for $t \downarrow 0$. If $0 < t < \bar{x}_j$, then $\bar{x}_j - t > 0$. This implies

$$u_i - \bar{x}_i = P_K(u)_i - P_K(\bar{x})_i = \begin{cases} 0, & \text{if } i \neq j, \\ -t, & \text{if } i = j, \end{cases} \text{ for } i = 1, 2, \ldots, n,$$

We have

$$\limsup_{u \to \bar{x}} \frac{\langle \theta, u - \bar{x}\rangle - \langle \bar{x}, P_K(u) - P_K(\bar{x})\rangle}{\|u - \bar{x}\| + \|P_K(u) - P_K(\bar{x})\|}$$

$$\geq \limsup_{t \downarrow 0, t < \bar{x}_j} \frac{-\langle \bar{x}, P_K(u) - P_K(\bar{x})\rangle}{\|u - \bar{x}\| + \|P_K(u) - P_K(\bar{x})\|}$$

$$= \limsup_{t \downarrow 0, t < \bar{x}_j} \frac{-(-t\bar{x}_j)}{2t}$$

$$= \frac{\bar{x}_j}{2} > 0.$$

It implies that

$$\bar{x}^+ \neq \emptyset \implies \theta \notin \widehat{D}^* P_K(\bar{x})(\bar{x}). \tag{4.12}$$

By (4.11) and (4.12), we proved that

$$\bar{x}^+ = \emptyset \iff \theta \in \widehat{D}^* P_K(\bar{x})(\bar{x}).$$

By (4.6), this implies (4.7). (4.8) immediately follows from (4.6) and (4.7). □

## 5. Mordukhovich derivatives of the metric projection operator in $l_2$

In this section, we consider the real Hilbert space $l_2$ with norm $\|\cdot\|$, with inner product $\langle \cdot, \cdot \rangle$ and with the origin $\theta$. Let $\mathbb{N}$ denote the set of all positive integers. To distinguish the differences between Euclidean spaces and the real Hilbert space $l_2$, we use some notations for the positive cone in real Hilbert space $l_2$ with different from the notations used in Euclidean spaces.

$$\mathbb{K} = \{x = (x_1, x_2, \ldots) \in l_2 : x_i \geq 0, \text{ for all } i \in \mathbb{N}\};$$

$$\mathbb{K}^+ = \{x = (x_1, x_2, \ldots) \in \mathbb{K}: x_i > 0, \text{ for all } i \in \mathbb{N}\};$$

$$\mathbb{K}^- = \{x = (x_1, x_2, \ldots) \in l_2: x_i < 0, \text{ for all } i \in \mathbb{N}\};$$

$$\widehat{\mathbb{K}} = \{x = (x_1, x_2, \ldots) \in l_2: |x_i| > 0, \text{ for all } i \in \mathbb{N} \text{ and}$$
$$\text{there are at least one pair } j, k \in \mathbb{N} \text{ with } x_j x_k < 0\}.$$

Let $N$ be a nonempty subset of $\mathbb{N}$ with complementary $\overline{N}$. We define some subsets in $l_2$ with respect to $N$.

$$\mathbb{R}^N = \{x = (x_1, x_2, \ldots) \in l_2: x_i = 0, \text{ for } i \in \overline{N}\};$$

$$\mathbb{K}_N = \{x = (x_1, x_2, \ldots) \in l_2: x_i \geq 0, \text{ for } i \in N\};$$

$$\partial \mathbb{K}_N = \{x = (x_1, x_2, \ldots) \in \mathbb{K}_N: x_i = 0, \text{ for } i \in N\};$$

$$\mathbb{Z}_N = \{x = (x_1, x_2, \ldots) \in l_2: x_i > 0, \text{ for all } i \in N \text{ and } x_i = 0, \text{ for all } i \in \overline{N}\}.$$

$\mathbb{R}^N$ is a closed subspace of $l_2$. $\mathbb{K}$ and $\mathbb{K}_N$ all are pointed closed and convex cones in $l_2$. We define an ordering relation $\leqslant_N$ on $l_2$, for any $y, z \in l_2$, by

$$z \leqslant_N y \quad \Longleftrightarrow \quad z_i \leq y_i, \text{ for all } i \in N \text{ and } y_i = z_i, \text{ for all } i \in \overline{N}.$$

We see that the ordering relation $\leqslant_N$ is a partial order on $l_2$. The interior of $\mathbb{K}$ is empty. It is proved in [8].

**Lemma 5.1 in [8]**. (a) *The interior of the positive cone $\mathbb{K}$ in $l_2$ is empty.*

(b) *If $N$ is an infinity set, then the interior of $\mathbb{K}_N$ is empty.*

Let $P_{\mathbb{K}}$ be the metric projection operator from $l_2$ onto $\mathbb{K}$. Similar, to Lemma 4.3 in [8], $P_{\mathbb{K}}$ has the following properties.

**Lemma 5.2 in [8]**. *Let $\mathbb{K}$ be the positive cone of $l_2$. $P_{\mathbb{K}}$ has the following properties.*

(a) *For any $x \in l_2$, $P_{\mathbb{K}}(x)$ is represented as follows*

$$P_{\mathbb{K}}(x)_i = \begin{cases} x_i, & \text{if } x_i > 0, \\ 0, & \text{if } x_i \leq 0, \end{cases} \text{ for } i \in \mathbb{N}. \tag{5.1}$$

(b) *$P_{\mathbb{K}}$ is positive homogeneous. For any $x \in l_2$,*

$$P_{\mathbb{K}}(\lambda x) = \lambda P_{\mathbb{K}}(x), \text{ for any } \lambda \geq 0. \tag{5.2}$$

Similarly, to (4.2) in section 4, for any given fixed $x \in \widehat{\mathbb{K}}$, we define a mapping $B(x; \cdot): l_2 \to l_2$, for any $w \in l_2$, by

$$(B(x; w))_i = \begin{cases} w_i, & \text{if } x_i > 0, \\ 0, & \text{if } x_i \leq 0, \end{cases} \text{ for all } i \in \mathbb{N}.$$

In this section, we study the generalized differentiability of the metric projection operator $P_\mathbb{K}: l_2 \to \mathbb{K}$. For the Fréchet differentiability of the metric projection operator $P_\mathbb{K}$, we recall the following results from [8].

**Theorem 5.3 in [8]**. *Let $\mathbb{K}$ be the positive cone of $l_2$. $P_\mathbb{K}$ has the following properties.*

(i) *In $\mathbb{K}^+$, we have*

  (a) *$P_\mathbb{K}$ is not Fréchet differentiability at any point in $\mathbb{K}^+$, that is,*
  $$\nabla P_K(x) \text{ does not exist, for any } x \in \mathbb{K}^+.$$

  (b) *$P_\mathbb{K}$ is Gâteaux directionally differentiable on $\mathbb{K}^+$ such that, for any $x \in \mathbb{K}^+$,*
  $$P'_K(x)(w) = w, \quad \text{for any } w \in l_2 \setminus \{\theta\}.$$

(ii) *In $\mathbb{K}^-$, we have*

  (a) *$P_\mathbb{K}$ is not Fréchet differentiability at any point in $\mathbb{K}^-$. that is,*
  $$\nabla P_K(x) \text{ does not exist, for any } x \in \mathbb{K}^-.$$

  (b) *$P_\mathbb{K}$ is Gâteaux directionally differentiable on $\mathbb{K}^-$ such that, for any $x \in \mathbb{K}^-$,*
  $$P'_K(x)(w) = \theta, \quad \text{for any } w \in \mathbb{R}^n \setminus \{\theta\}.$$

(iii) *In $\widehat{\mathbb{K}}$, we have*

  (a) *$P_\mathbb{K}$ is not Fréchet differentiable on $\widehat{\mathbb{K}}$, that is,*
  $$\nabla P_\mathbb{K}(x) \text{ does not exist, for any } x \in \widehat{\mathbb{K}}.$$

  (b) *$P_\mathbb{K}$ is Gâteaux directionally differentiable on $\widehat{\mathbb{K}}$ such that, for any $x \in \widehat{\mathbb{K}}$,*
  $$P'_\mathbb{K}(x)(w) = B(x; w), \quad \text{for any } w \in l_2 \setminus \{\theta\}.$$

**Theorem 5.1**. *Let $M$ be a nonempty finite subset of $\mathbb{N}$ with complementary $\bar{M}$. The Mordukhovich derivatives of $P_\mathbb{K}$ have the following representations.*

(i) *For any $\bar{x} \in l_2$, we have*
$$\widehat{D}^* P_\mathbb{K}(\bar{x})(\theta) = \{\theta\};$$

(ii) *Let $\bar{x} \in \mathbb{Z}_M$. For any $y \in l_2$,*
$$y \in \mathbb{K}_{\bar{M}} \iff y \in \widehat{D}^* P_\mathbb{K}(\bar{x})(y). \tag{5.3}$$

(iii) *Let $\bar{x} \in \mathbb{Z}_M$. For any $y \in \mathbb{K}_{\bar{M}}$, we have*
$$\widehat{D}^* P_\mathbb{K}(\bar{x})(y) = \{z \in \mathbb{K}_{\bar{M}}: z \preccurlyeq_{\bar{M}} y\}. \tag{5.4}$$

*In particular,*

$$y \in \partial \mathbb{K}_{\bar{M}} \implies \widehat{D}^* P_{\mathbb{K}}(\bar{x})(y) = \{y\}.$$

*Proof.* Proof of (i). The proof of (i) is same with the proof of part (a) of (iv) in Theorem 4.1. So, the proof of (i) is omitted here.

We prove (ii). For any $\bar{x} \in \mathbb{Z}_M$, we first prove

$$y \in \mathbb{K}_{\bar{M}} \implies y \in \widehat{D}^* P_{\mathbb{K}}(\bar{x})(y). \tag{5.3.1}$$

Let $\delta = \min\{\bar{x}_i : \text{for all } i \in M\}$. By $\bar{x} \in \mathbb{Z}_M$, we have $\delta > 0$. Then, for any $u = (u_1, u_2, \ldots) \in l_2$,

$$|u - \bar{x}| < \frac{\delta}{2} \implies u_i > 0, \text{ for all } i \in M.$$

This implies that, if $|u - \bar{x}| < \frac{\delta}{2}$, then

$$P_{\mathbb{K}}(u)_i = u_i, \text{ for all } i \in M,$$

$$P_{\mathbb{K}}(\bar{x})_i = \bar{x}_i, \text{ for all } i \in M,$$

and
$$\bar{x}_i = 0 \text{ and } P_{\mathbb{K}}(\bar{x})_i = 0, \text{ for all } i \in \bar{M}. \tag{5.5}$$

Then, for any $u = (u_1, u_2, \ldots) \in l_2$ satisfying $|u - \bar{x}| < \frac{\delta}{2}$, for $y \in \mathbb{K}_{\bar{M}}$, by (5.5), we have

$$\langle y, u - \bar{x} \rangle - \langle y, P_{\mathbb{K}}(u) - P_{\mathbb{K}}(\bar{x}) \rangle$$

$$= \langle y, u - P_{\mathbb{K}}(u) \rangle - \langle y, \bar{x} - P_{\mathbb{K}}(\bar{x}) \rangle$$

$$= \Sigma_{i \in M} y_i (u_i - P_{\mathbb{K}}(u)_i) + \Sigma_{i \in \bar{M}} y_i (u_i - P_{\mathbb{K}}(u)_i) - \Sigma_{i \in M} y_i (\bar{x}_i - P_{\mathbb{K}}(\bar{x})_i) - \Sigma_{i \in \bar{M}} y_i (\bar{x}_i - P_{\mathbb{K}}(\bar{x})_i)$$

$$= \Sigma_{i \in M} y_i 0 + \Sigma_{i \in \bar{M}} y_i (u_i - P_{\mathbb{K}}(u)_i) - \Sigma_{i \in M} y_i 0 - \Sigma_{i \in \bar{M}} y_i (0 - 0)$$

$$= \Sigma_{i \in \bar{M}} y_i (u_i - P_{\mathbb{K}}(u)_i)$$

$$= \Sigma_{i \in \bar{M}, u_i \leq 0} y_i (u_i - P_{\mathbb{K}}(u)_i) + \Sigma_{i \in \bar{M}, u_i > 0} y_i (u_i - P_{\mathbb{K}}(u)_i)$$

$$= \Sigma_{i \in \bar{M}, u_i \leq 0} y_i (u_i - 0) + \Sigma_{i \in \bar{M}, u_i > 0} y_i 0$$

$$= \Sigma_{i \in \bar{M}, u_i \leq 0} y_i u_i$$

$$\leq 0 \text{ (It is because that } y_i \geq 0, \bar{x}_i = 0, \text{ for all } i \in \bar{M}). \tag{5.6}$$

It follows that

$$\limsup_{u \to \bar{x}} \frac{\langle y, u - \bar{x} \rangle - \langle y, P_{\mathbb{K}}(u) - P_{\mathbb{K}}(\bar{x}) \rangle}{\|u - \bar{x}\| + \|P_{\mathbb{K}}(u) - P_{\mathbb{K}}(\bar{x})\|}$$

$$= \limsup_{u \to \bar{x}, |u-\bar{x}| < \frac{\delta}{2}} \frac{\langle y, u-\bar{x} \rangle - \langle y, P_{\mathbb{K}}(u) - P_{\mathbb{K}}(\bar{x}) \rangle}{\|u-\bar{x}\| + \|P_{\mathbb{K}}(u) - P_{\mathbb{K}}(\bar{x})\|}$$

$$= \limsup_{u \to \bar{x}, |u-\bar{x}| < \frac{\delta}{2}} \frac{\sum_{i \in \bar{M}, u_i \le 0} y_i u_i}{\|u-\bar{x}\| + \|P_{\mathbb{K}}(u) - P_{\mathbb{K}}(\bar{x})\|}$$

$$\le 0.$$

This implies that, for any $\bar{x} \in \mathbb{Z}_M$, we have

$$y \in \mathbb{K}_{\bar{M}} \implies y \in \widehat{D}^* P_{\mathbb{K}}(\bar{x})(y). \tag{5.3.1}$$

For fixed $\bar{x} \in \mathbb{Z}_M$ and for $y \in l_2$, suppose $y \notin \mathbb{K}_{\bar{M}}$. Then, there is a positive integer $m$ with $m \in \bar{M}$ such that $y_m < 0$. We take a directional line segment $u$ to approach to $\bar{x}$ as follows:

$$u_i = \begin{cases} \bar{x}_i, & \text{for } i \in M, \\ -t, & \text{for } i = m, \\ 0, & \text{for } i \in \bar{M} \setminus \{m\}. \end{cases}$$

For $t > 0$, by $\bar{x}_m = 0$ (It is because that $\bar{x} \in \mathbb{Z}_M$ and $m \in \bar{M}$), one has

$$P_{\mathbb{K}}(u) = P_{\mathbb{K}}(\bar{x}) \text{ and } \langle y, u - \bar{x} \rangle = -t y_m > 0.$$

We calculate

$$\limsup_{u \to \bar{x}} \frac{\langle y, u-\bar{x} \rangle - \langle y, P_{\mathbb{K}}(u) - P_{\mathbb{K}}(\bar{x}) \rangle}{\|u-\bar{x}\| + \|P_{\mathbb{K}}(u) - P_{\mathbb{K}}(\bar{x})\|}$$

$$\ge \limsup_{t \downarrow 0} \frac{\langle y, u-\bar{x} \rangle - \langle y, \theta \rangle}{\|u-\bar{x}\| + \|\theta\|}$$

$$= \limsup_{t \downarrow 0} \frac{-t y_m}{t}$$

$$= -y_m > 0.$$

This proves that, for any $\bar{x} \in \mathbb{Z}_M$, we have

$$y \notin \mathbb{K}_{\bar{M}} \implies y \notin \widehat{D}^* P_{\mathbb{K}}(\bar{x})(y). \tag{5.3.2}$$

Hence, (5.3) is proved by (5.3.1) and (5.3.2), which proves (ii).

Proof of (iii). For any fixed $\bar{x} \in \mathbb{Z}_M$ and for any $y, z \in l_2$, suppose that $z \in \mathbb{K}_{\bar{M}}$ and $z \preceq_{\bar{M}} y$. We prove

$$z \in \widehat{D}^* P_{\mathbb{K}}(\bar{x})(y). \tag{5.4.1}$$

Since $z \in \mathbb{K}_{\bar{M}}$ and $z \preceq_{\bar{M}} y$, then $y \in \mathbb{K}_{\bar{M}}$. By (ii), it is proved that $y \in \widehat{D}^* P_{\mathbb{K}}(\bar{x})(y)$.

Since it is assumed that $z \in \mathbb{K}_{\bar{M}}$, by (ii), we have $z \in \widehat{D}^* P_{\mathbb{K}}(\bar{x})(z)$. This implies that

$$\limsup_{u \to \bar{x}} \frac{\langle z, u - \bar{x} \rangle - \langle z, P_{\mathbb{K}}(u) - P_{\mathbb{K}}(\bar{x}) \rangle}{\|u - \bar{x}\| + \|P_{\mathbb{K}}(u) - P_{\mathbb{K}}(\bar{x})\|} \leq 0. \tag{5.4.2}$$

Let $\delta = \min\{\bar{x}_i : \text{for all } i \in M\}$ as defined in the proof of (ii). Then, for any $u = (u_1, u_2, \ldots) \in l_2$,

$$|u - \bar{x}| < \frac{\delta}{2} \quad \Longrightarrow \quad u_i > 0, \text{ for all } i \in M.$$

This implies that, if $|u - \bar{x}| < \frac{\delta}{2}$, then

$$P_{\mathbb{K}}(u)_i = u_i > 0, \text{ for all } i \in M,$$

$$P_{\mathbb{K}}(\bar{x})_i = \bar{x}_i > 0, \text{ for all } i \in M,$$

and
$$\bar{x}_i = 0 \text{ and } P_{\mathbb{K}}(\bar{x})_i = 0, \text{ for all } i \in \bar{M}. \tag{5.5}$$

Then, for any $u = (u_1, u_2, \ldots) \in l_2$ satisfying $|u - \bar{x}| < \frac{\delta}{2}$, for $z \in \mathbb{K}_{\bar{M}}$ with $z \leqslant_{\bar{M}} y$, similar, to the calculation of (5.6), this implies

$$\langle z, u - \bar{x} \rangle - \langle y, P_{\mathbb{K}}(u) - P_{\mathbb{K}}(\bar{x}) \rangle$$

$$= \langle z, u - \bar{x} \rangle - \langle z, P_{\mathbb{K}}(u) - P_{\mathbb{K}}(\bar{x}) \rangle + \langle z, P_{\mathbb{K}}(u) - P_{\mathbb{K}}(\bar{x}) \rangle - \langle y, P_{\mathbb{K}}(u) - P_{\mathbb{K}}(\bar{x}) \rangle$$

$$= \langle z, u - \bar{x} \rangle - \langle z, P_{\mathbb{K}}(u) - P_{\mathbb{K}}(\bar{x}) \rangle + \langle z - y, P_{\mathbb{K}}(u) - P_{\mathbb{K}}(\bar{x}) \rangle$$

$$= \langle z, u - \bar{x} \rangle - \langle z, P_{\mathbb{K}}(u) - P_{\mathbb{K}}(\bar{x}) \rangle + \sum_{i \in M}(z_i - y_i)(P_{\mathbb{K}}(u)_i - P_{\mathbb{K}}(\bar{x})_i) + \sum_{i \in \bar{M}}(z_i - y_i)(P_{\mathbb{K}}(u)_i - P_{\mathbb{K}}(\bar{x})_i)$$

$$= \langle z, u - \bar{x} \rangle - \langle z, P_{\mathbb{K}}(u) - P_{\mathbb{K}}(\bar{x}) \rangle + \sum_{i \in M} 0 (P_{\mathbb{K}}(u)_i - P_{\mathbb{K}}(\bar{x})_i) + \sum_{i \in \bar{M}}(z_i - y_i)(P_{\mathbb{K}}(u)_i - 0)$$

$$= \langle z, u - \bar{x} \rangle - \langle z, P_{\mathbb{K}}(u) - P_{\mathbb{K}}(\bar{x}) \rangle + \sum_{i \in \bar{M}}(z_i - y_i) P_{\mathbb{K}}(u)_i$$

$$= \langle z, u - \bar{x} \rangle - \langle z, P_{\mathbb{K}}(u) - P_{\mathbb{K}}(\bar{x}) \rangle + \sum_{i \in \bar{M}, u_i \leq 0}(z_i - y_i) P_{\mathbb{K}}(u)_i + \sum_{i \in \bar{M}, u_i > 0}(z_i - y_i) P_{\mathbb{K}}(u)_i$$

$$= \langle z, u - \bar{x} \rangle - \langle z, P_{\mathbb{K}}(u) - P_{\mathbb{K}}(\bar{x}) \rangle + \sum_{i \in \bar{M}, u_i \leq 0}(z_i - y_i) 0 + \sum_{i \in \bar{M}, u_i > 0}(z_i - y_i) u_i$$

$$= \langle z, u - \bar{x} \rangle - \langle z, P_{\mathbb{K}}(u) - P_{\mathbb{K}}(\bar{x}) \rangle + \sum_{i \in \bar{M}, u_i > 0}(z_i - y_i) u_i.$$

Since $z_i - y_i \leq 0$, for all $i \in \bar{M}$, this implies

$$\sum_{i \in \bar{M}, u_i > 0}(z_i - y_i) u_i \leq 0. \tag{5.4.3}$$

By (5.4.2) and (5.4.3), it follows that

$$\limsup_{u \to \bar{x}} \frac{\langle z, u - \bar{x} \rangle - \langle y, P_{\mathbb{K}}(u) - P_{\mathbb{K}}(\bar{x}) \rangle}{\|u - \bar{x}\| + \|P_{\mathbb{K}}(u) - P_{\mathbb{K}}(\bar{x})\|}$$

$$\begin{aligned}
&= \limsup_{u \to \bar{x}, |u-\bar{x}| < \frac{\delta}{2}} \frac{\langle z, u-\bar{x} \rangle - \langle y, P_{\mathbb{K}}(u) - P_{\mathbb{K}}(\bar{x}) \rangle}{\|u-\bar{x}\| + \|P_{\mathbb{K}}(u) - P_{\mathbb{K}}(\bar{x})\|} \\
&= \limsup_{u \to \bar{x}, |u-\bar{x}| < \frac{\delta}{2}} \frac{\langle z, u-\bar{x} \rangle - \langle z, P_{\mathbb{K}}(u) - P_{\mathbb{K}}(\bar{x}) \rangle + \sum_{i \in \bar{M}, u_i > 0}(z_i - y_i)u_i}{\|u-\bar{x}\| + \|P_{\mathbb{K}}(u) - P_{\mathbb{K}}(\bar{x})\|} \\
&\leq \limsup_{u \to \bar{x}, |u-\bar{x}| < \frac{\delta}{2}} \frac{\langle z, u-\bar{x} \rangle - \langle z, P_{\mathbb{K}}(u) - P_{\mathbb{K}}(\bar{x}) \rangle}{\|u-\bar{x}\| + \|P_{\mathbb{K}}(u) - P_{\mathbb{K}}(\bar{x})\|} + \limsup_{u \to \bar{x}, |u-\bar{x}| < \frac{\delta}{2}} \frac{\sum_{i \in \bar{M}, u_i > 0}(z_i - y_i)u_i}{\|u-\bar{x}\| + \|P_{\mathbb{K}}(u) - P_{\mathbb{K}}(\bar{x})\|} \\
&\leq 0 + \limsup_{u \to \bar{x}, |u-\bar{x}| < \frac{\delta}{2}} \frac{0}{\|u-\bar{x}\| + \|P_{\mathbb{K}}(u) - P_{\mathbb{K}}(\bar{x})\|} \\
&= 0.
\end{aligned}$$

This implies that, for fixed $\bar{x} \in \mathbb{Z}_M$ and $y, z \in l_2$, if $z \in \mathbb{K}_{\bar{M}}$ and $z \preccurlyeq_{\bar{M}} y$, then, $z \in \widehat{D}^* P_{\mathbb{K}}(\bar{x})(y)$, which proves (5.4.1). That is,

For fixed $\bar{x} \in \mathbb{Z}_M$ and for any $y, z \in l_2$ with $y \in \mathbb{K}_{\bar{M}}$. Suppose $z$ does not satisfy the conditions that $z \in \mathbb{K}_{\bar{M}}$ and $z \preccurlyeq_{\bar{M}} y$. This is equivalent to that $z$ satisfies the following conditions:

$$\text{either, } z \in \mathbb{K}_{\bar{M}} \text{ and } z \not\preccurlyeq_{\bar{M}} y, \quad \text{or,} \quad z \notin \mathbb{K}_{\bar{M}}. \tag{5.7}$$

Conditions in (5.7) are equivalent to: there are numbers $j$, $m$, $n$ with $j \in M$ and $m, n \in \bar{M}$, such that at least one of the following conditions is satisfied:

$$z_j \neq y_j, \text{ or } z_n > y_n, \text{ or } z_m < 0. \tag{5.8}$$

For $z \in l_2$, we prove case by case that if $z$ satisfies any one of the 3 conditions in (5.8), then we have $z \notin \widehat{D}^* P_{\mathbb{K}}(\bar{x})(y)$.

Case 1. Assume that there is $j \in M$ satisfying $z_j \neq y_j$. Let $\delta = \min\{\bar{x}_i: \text{ for all } i \in M\}$ as defined in the proof of (ii). We take a directional line segment $u$ to approach to $\bar{x}$ as follows:

$$u_i = \begin{cases} \bar{x}_i, & \text{for } i \in M \setminus \{j\}, \\ \begin{cases} \bar{x}_j - t, \text{ if } y_j > z_j \\ \bar{x}_j + t, \text{ if } y_j < z_j \end{cases}, & \text{for } i = j, \\ 0, & \text{for } i \in \bar{M}. \end{cases} \tag{5.9}$$

For $0 < t < \frac{1}{2}\delta$, we have

$$\langle z, u - \bar{x} \rangle - \langle y, P_{\mathbb{K}}(u) - P_{\mathbb{K}}(\bar{x}) \rangle = \begin{cases} t(y_j - z_j), & \text{if } y_j > z_j \\ t(z_j - y_j), & \text{if } y_j < z_j \end{cases}. \tag{5.10}$$

Then, for any $u = (u_1, u_2, \ldots) \in l_2$ defined by (5.9), by (5.10), we calculate

$$\operatorname*{limsup}_{u\to\bar{x}} \frac{\langle z,u-\bar{x}\rangle - \langle y, P_{\mathbb{K}}(u) - P_{\mathbb{K}}(\bar{x})\rangle}{\|u-\bar{x}\| + \|P_{\mathbb{K}}(u) - P_{\mathbb{K}}(\bar{x})\|}$$

$$\geq \operatorname*{limsup}_{t\downarrow 0} \frac{\langle z,u-\bar{x}\rangle - \langle y, P_{\mathbb{K}}(u) - P_{\mathbb{K}}(\bar{x})\rangle}{\|u-\bar{x}\| + \|P_{\mathbb{K}}(u) - P_{\mathbb{K}}(\bar{x})\|}$$

$$= \operatorname*{limsup}_{t\downarrow 0} \frac{t|y_j - z_j|}{2t}$$

$$= \frac{|y_j - z_j|}{2}$$

$$> 0.$$

This proves that, for any $\bar{x} \in \mathbb{Z}_M$, if there is $j \in M$ satisfying $z_j \neq y_j$, then

$$z \notin \widehat{D}^* P_{\mathbb{K}}(\bar{x})(y). \tag{5.8.1}$$

Case 2. Assume that there is $n \in \bar{M}$ satisfying $z_n > y_n$ (here $y_n \geq 0$, for $n \in \bar{M}$). We take a directional line segment $u$ to approach to $\bar{x}$ as follows:

$$u_i = \begin{cases} \bar{x}_i, & \text{for } i \in M, \\ t, & \text{for } i = n, \\ 0, & \text{for } i \in \bar{M}\setminus\{n\}. \end{cases} \tag{5.11}$$

For $t > 0$, by $\bar{x}_n = 0$, for $n \in \bar{M}$, we have

$$\langle z, u - \bar{x}\rangle - \langle y, P_{\mathbb{K}}(u) - P_{\mathbb{K}}(\bar{x})\rangle$$

$$= z_n(t - 0) - t y_n$$

$$= t(z_n - y_n). \tag{5.12}$$

Then, for any $u = (u_1, u_2, \ldots) \in l_2$ defined by (5.11), by (5.12), by $\bar{x}_n = 0$, $z_n > y_n \geq 0$, we calculate

$$\operatorname*{limsup}_{u\to\bar{x}} \frac{\langle y,u-\bar{x}\rangle - \langle y, P_{\mathbb{K}}(u) - P_{\mathbb{K}}(\bar{x})\rangle}{\|u-\bar{x}\| + \|P_{\mathbb{K}}(u) - P_{\mathbb{K}}(\bar{x})\|}$$

$$\geq \operatorname*{limsup}_{t\downarrow 0} \frac{\langle y,u-\bar{x}\rangle - \langle y, P_{\mathbb{K}}(u) - P_{\mathbb{K}}(\bar{x})\rangle}{\|u-\bar{x}\| + \|P_{\mathbb{K}}(u) - P_{\mathbb{K}}(\bar{x})\|}$$

$$= \operatorname*{limsup}_{t\downarrow 0} \frac{z_n(t-0) - t y_n}{z_n(t - \bar{x}_n) + t}$$

$$= \operatorname*{limsup}_{t\downarrow 0} \frac{t(z_n - y_n)}{t(1 + z_n) - z_n \bar{x}_n}$$

$$= z_n - y_n$$

$$> 0 \text{ (It is because that } \bar{x}_n = 0, z_n > y_n \geq 0\text{)}.$$

This proves that, for any $\bar{x} \in \mathbb{Z}_M$, if there is $n \in \bar{M}$ satisfying $z_n > y_n$, then

$$z \notin \widehat{D}^* P_\mathbb{K}(\bar{x})(y). \tag{5.8.2}$$

Case 3. Assume that there is $m \in \bar{M}$ satisfying $z_m < 0$. We take a directional line segment $u$ to approach to $\bar{x}$ as follows:

$$u_i = \begin{cases} \bar{x}_i, & \text{for } i \in M, \\ -t, & \text{for } i = m \\ 0, & \text{for } i \in \bar{M}\setminus\{m\}, \end{cases}, \text{ for } i = 1, 2, \dots. \tag{5.13}$$

For $t > 0$, we have

$$\langle z, u - \bar{x}\rangle - \langle y, P_\mathbb{K}(u) - P_\mathbb{K}(\bar{x})\rangle$$

$$= \langle z, u - \bar{x}\rangle - \langle y, \theta\rangle$$

$$= -tz_m. \tag{5.14}$$

Then, for any $u = (u_1, u_2, \dots) \in l_2$ defined by (5.13), by (5.14), we calculate

$$\limsup_{u \to \bar{x}} \frac{\langle y, u-\bar{x}\rangle - \langle y, P_\mathbb{K}(u) - P_\mathbb{K}(\bar{x})\rangle}{\|u-\bar{x}\| + \|P_\mathbb{K}(u) - P_\mathbb{K}(\bar{x})\|}$$

$$\geq \limsup_{t \downarrow 0} \frac{\langle y, u-\bar{x}\rangle - \langle y, P_\mathbb{K}(u) - P_\mathbb{K}(\bar{x})\rangle}{\|u-\bar{x}\| + \|\theta\|}$$

$$= \limsup_{t \downarrow 0} \frac{-tz_m}{t}$$

$$= -z_m$$

$$> 0.$$

This proves that, for any $\bar{x} \in \mathbb{Z}_M$, if there is $m \in \bar{M}$ satisfying $z_m < 0$, then

$$z \notin \widehat{D}^* P_\mathbb{K}(\bar{x})(y). \tag{5.8.3}$$

By (5.8.1), (5.8.2) and (5.8.3), we proved that for any $\bar{x} \in \mathbb{Z}_M$, $z \in l_2$ and $y \in \mathbb{K}_{\bar{M}}$, if $z \not\leq_{\bar{M}} y$, then,

$$z \notin \widehat{D}^* P_\mathbb{K}(\bar{x})(y). \tag{5.8.4}$$

Then, (5.4) follows from (5.4.1) and (5.8.4). This proves (iii) of this theorem.

Finally, we prove the special case in (iii):

$$y \in \partial \mathbb{K}_{\bar{M}} \implies \widehat{D}^* P_\mathbb{K}(\bar{x})(y) = \{y\}.$$

We proved that for $\bar{x} \in \mathbb{Z}_M$ and for any $y \in \mathbb{K}_{\bar{M}}$, we have

$$\widehat{D}^* P_\mathbb{K}(\bar{x})(y) = \{z \in \mathbb{K}_{\bar{M}}: z \leq_{\bar{M}} y\}. \tag{5.4}$$

By definition,

$$y \in \partial \mathbb{K}_{\bar{M}} \implies y_i = 0, \text{ for all } i \notin M.$$

For any $z \in \mathbb{K}_{\bar{M}}$ with $z \preccurlyeq_{\bar{M}} y$, we have

$$z_i = y_i, \text{ for all } i \in M \text{ and } 0 \leq z_i \leq y_i = 0, \text{ for all } i \notin M.$$

This implies $z = y$; and therefore, in (5.4), $\{z \in \mathbb{K}_{\bar{M}} : z \preccurlyeq_{\bar{M}} y\} = \{y\}$. □

**Corollary 5.2.** *Let $M$ be a nonempty finite subset of $\mathbb{N}$ with complementary $\bar{M}$. Then, $P_{\mathbb{K}}$ is not Fréchet differentiability on $\mathbb{Z}_M$. That is, for any $\bar{x} \in \mathbb{Z}_M$, $\nabla P_{\mathbb{K}}(\bar{x})$ does not exist.*

*Proof.* For an arbitrary $\bar{x} \in \mathbb{Z}_M$, by part (iii) in Theorem 5.1, we have

$$\widehat{D}^* P_{\mathbb{K}}(\bar{x})(y) = \{z \in \mathbb{K}_{\bar{M}} : z \preccurlyeq_{\bar{M}} y\}, \text{ for any } y \in \mathbb{K}_{\bar{M}}.$$

This set is not a singleton for $y \in \mathbb{K}_{\bar{M}}$, in general. It is a singleton, only if $y \in \partial \mathbb{K}_{\bar{M}}$, in which $\widehat{D}^* P_{\mathbb{K}}(\bar{x})(y) = \{y\}$, which is a singleton.

If $P_{\mathbb{K}}$ is Fréchet differentiability at $\bar{x} \in \mathbb{Z}_M$ with Fréchet derivative $\nabla P_{\mathbb{K}}(\bar{x})$, then, by Theorem 1.38 in [11], we should have

$$\widehat{D}^* P_{r\mathbb{B}}(\bar{x})(y) = \{\nabla P_{r\mathbb{B}}(\bar{x})(y)\},$$

That is always a singleton for any $y \in \mathbb{K}_{\bar{M}}$. This contradiction shows that $P_{\mathbb{K}}$ is not Fréchet differentiability at $\bar{x} \in \mathbb{Z}_M$. □

Let $M$ be a nonempty finite subset of $\mathbb{N}$ with complementary $\bar{M}$. Recall the notations:

$$\mathbb{R}^M = \{x = (x_1, x_2, \ldots) \in l_2 : x_i = 0, \text{ for all } i \notin M\};$$

$$\mathbb{Z}_M = \{x = (x_1, x_2, \ldots) \in l_2 : x_i > 0, \text{ for all } i \in M \text{ and } x_i = 0, \text{ for all } i \notin M\}.$$

In particular, if $N = \{1, 2, \ldots, n\}$, for some $n \in \mathbb{N}$, then $\mathbb{R}^N$ is abbreviated as $\mathbb{R}^n$.

For the given nonempty finite subset $M$ of $\mathbb{N}$, $\mathbb{R}^M$ is a proper closed subspace of $l_2$. Let $I_{\mathbb{R}^M}$ denote the identity mapping on $\mathbb{R}^M$. Next, we investigate the Mordukhovich derivatives of the metric projection operator $P_{\mathbb{K}} : l_2 \to \mathbb{K}$ on the subspace $\mathbb{R}^M$ of $l_2$.

**Corollary 5.3.** *Let $M$ be a nonempty finite subset of $\mathbb{N}$ with complementary $\bar{M}$. We have*

$$\widehat{D}^* P_{\mathbb{K}} \big|_{\mathbb{Z}_M} = I_{\mathbb{R}^M}.$$

*That is, for any $\bar{x} \in \mathbb{Z}_M$, we have*

$$\widehat{D}^* P_{\mathbb{K}}(\bar{x})(y) = y, \text{ for any } y \in \mathbb{R}^M. \tag{5.15}$$

*Proof.* Notice that $\partial \mathbb{K}_{\bar{M}} = \mathbb{R}^M$ and $\mathbb{Z}_M \subseteq \mathbb{R}^M$. Let $\bar{x} \in \mathbb{Z}_M$ be arbitrarily given. For any $y \in \mathbb{R}^M$,

since $\partial \mathbb{K}_{\bar{M}} = \mathbb{R}^M$, by part (iii) in Theorem 5.1, we have

$$y \in \mathbb{R}^M = \partial \mathbb{K}_{\bar{M}} \implies \widehat{D}^* P_\mathbb{K}(\bar{x})(y) = \{y\} = y.$$

This proves (5.15). □

### 6. Conclusion and remarks

In section 3 of this paper, we study the generalized differentiability of the metric projection onto closed balls centered at the origin in Hilbert spaces. In [6, 7], the directional differentiability of the metric projection onto closed balls is studied in uniformly convex and uniformly smooth Banach spaces and Hilbert spaces, in which the considered balls have center at arbitrarily given point $c$ in the spaces.

Similar to Theorem 4.2 in [7] recalled in Section 2, we think that the results about Mordukhovich derivatives of $P_{r\mathbb{B}}$ proved in Theorem 3.2 in this paper can be also extended to metric projection operator onto any closed balls $\mathbb{B}(c, r)$, which has center $c$ in $H$ and with radius $r > 0$.

### Acknowledgments

The author is very grateful to Professor Boris S. Mordukhovich for his kind communications, valuable suggestions and enthusiasm encouragements in the development stage of this paper.### References